\documentclass[a4paper, reqno]{amsart}
\usepackage{graphicx} % Required for inserting images

\title[Poset functor cocalculus and applications to TDA]{Poset functor cocalculus and applications to topological data analysis}
\author{Bjørnar Gullikstad Hem}

\usepackage{amsmath}
\usepackage{amsthm}
\usepackage{amssymb}
\usepackage[hidelinks,colorlinks=true]{hyperref}
\hypersetup{allcolors=[rgb]{0,0.31,0.62}}
\usepackage{thmtools}   % Lar deg bruke \autoref.
\usepackage{mathrsfs}
\usepackage{enumitem}
\usepackage{mathtools}
\usepackage{caption}
\usepackage{subcaption}

\usepackage{tikz}
\usetikzlibrary{matrix} % For å kunne bruke \diagram må du ha denne.
\newcommand{\diagram}[3]{\matrix (#1) [matrix of math nodes,row
  sep={#2},column sep={#3},text height=1.5ex,text
  depth=0.25ex]}
\usepackage{tikz-cd}

\theoremstyle{plain}
\newtheorem{theorem}{Theorem}[section]
\newtheorem{proposition}[theorem]{Proposition}
\newtheorem{corollary}[theorem]{Corollary}
\newtheorem{lemma}[theorem]{Lemma}

\theoremstyle{definition}
\newtheorem{definition}[theorem]{Definition}
\newtheorem{example}[theorem]{Example}
\newtheorem{remark}[theorem]{Remark}

\DeclareMathOperator{\id}{id}
\DeclareMathOperator{\colim}{colim}
\DeclareMathOperator{\Hom}{Hom}
\newcommand{\Hhom}{\mathbf{Hom}}
\DeclareMathOperator{\op}{op}

\newcommand{\R}{\mathbb{R}} % Tavlefet R - \R
\newcommand{\Z}{\mathbb{Z}}

\newcommand{\Fb}{\mathbb{Fb}}

\newcommand{\C}{\mathscr{C}}
\newcommand{\M}{\mathscr{M}}

\newcommand{\bigo}{\mathcal{O}}
\newcommand{\SCplx}{\mathsf{SCplx}}
\newcommand{\Spaces}{\mathsf{Spaces}}
\newcommand{\Vecs}{\mathsf{Vec}}

\newcommand{\Po}{\mathcal{P}}
\newcommand{\X}{\mathcal{X}}
\newcommand{\Yc}{\mathcal{Y}}
\newcommand{\Zc}{\mathcal{Z}}
\newcommand{\Ch}{\mathrm{Ch}}
\newcommand{\Nn}{\mathbb{N}}
\newcommand{\nnR}{\mathbb{R}_{\ge 0}} % nonnegative reals
\newcommand{\Trans}{\mathbf{Trans}}
\newcommand{\Ib}{\mathbf{I}}
\newcommand{\Lke}{\mathbf{L}}
\newcommand{\Lv}{\Lambda_{\varepsilon}}

\DeclareMathOperator{\holim}{holim}
\DeclareMathOperator{\hocolim}{hocolim}

\DeclareMathOperator{\Fun}{Fun}
\DeclareMathOperator{\const}{const}
\DeclareMathOperator{\Ho}{Ho}
\DeclareMathOperator{\Lan}{Lan}
\DeclareMathOperator{\hoLan}{hoLan}
\DeclareMathOperator{\jdim}{jdim}

\numberwithin{equation}{subsection}

\newtheorem{thmx}{Theorem}

\newtheorem*{proposition*}{Proposition}

\begin{document}

\maketitle

\begin{abstract}
%Functor calculus is a general framework for studying functors, and comes in many flavors. Common for all of them is that one defines a notion of degree $n$ functors, and then studies a functor by considering its degree $n$ approximations. In Goodwillie’s functor calculus, the degree $n$ functors are $n$--excisive functors, i.e., functors that take strongly cocartesian $(n+1)$--cubes to cartesian $(n+1)$--cubes. 
%Inspired by Weiss’ manifold calculus, we construct a flavor of functor calculus for studying functors from a distributive lattice to a model category. In this calculus, the degree $n$ functors take strongly bicartesian $(n+1)$--cubes to cocartesian $(n+1)$--cubes.
%Our motivation is to better understand multi-persistence modules, as these modules can be viewed as functors from $\Nn^k$ into a category of chain complexes. It is a well-known problem in TDA that multi-persistence modules do not always admit an interval decomposition, and hence data from multi-persistence homology is difficult to analyze.
%We show that in our calculus, each functor out of $\Nn^k$ is degree $k$ and has a well-behaved telescope of lower degree approximations.
We introduce a new flavor of functor cocalculus, called \emph{poset cocalculus}, as a tool for studying approximations in topological data analysis. Given a functor from a distributive lattice to a model category, poset cocalculus produces a Taylor telescope of codegree $n$ approximations of the functor, where a codegree $n$ functor takes strongly bicartesian $(n+1)$--cubes to homotopy cocartesian $(n+1)$--cubes.
We give several applications of this new functor cocalculus.
We prove that the codegree $n$ approximation of a multipersistence module is stable under an appropriate notion of interleaving distance.
We draw connections to filtrations of simplicial complexes, and show that the Vietoris-Rips filtration is precisely the codegree 2 approximation of the Čech filtration.
We demonstrate that the codegree 1 approximation of the space of simplicial maps between two simplicial complexes is in some sense the space of continuous maps between their realizations, and that this statement can be made precise.
\end{abstract}

\setcounter{tocdepth}{1}
\tableofcontents

\section{Introduction}

\subsection{Background}

\subsubsection{Topological data analysis}

Topological data analysis (TDA) can be described as the general approach to applying methods from topology to data analysis.
One of the most common methods in TDA is \emph{persistent homology}. Given a point cloud, persistent homology provides successive approximations to the homology groups of the underlying space, by constructing an $\R_{\ge 0}$--parameterized filtration of simplicial complexes from the point cloud, and then post-composing with the functor $H_i( - ; \Fb)$, i.e., homology in some degree $i$ over some field $\Fb$. In total, this gives a functor
\begin{equation*}
    F \colon \nnR \to \Vecs_{\Fb}
\end{equation*}
from the poset $\R_{\ge 0}$ (i.e., nonnegative real numbers) to $\Vecs_{\Fb}$ (vector spaces over $\Fb$), which we call a \emph{persistence module}.

There are several methods of constructing an $\nnR$-parametrized filtration from a point cloud. Two of the most common ones are the \emph{Čech filtration} and its less computationally expensive approximation, the \emph{Vietoris-Rips filtration}. Given a point cloud $V \subset \R^d$, the \emph{Čech filtration} at $t \in \nnR$ is the simplicial complex $\{\sigma \subset V : f_{\mathbf{C}}(\sigma) \le t\}$, where $f_{\mathbf{C}}$ is the function that sends a simplex to its radius (i.e., the radius of the smallest ball containing the simplex). The Vietoris-Rips complex is defined similarly, except that $f_{\mathbf{C}}$ is replaced by $f_{\mathbf{VR}}$, which sends a simplex to its diameter (i.e., the maximum distance between two points in the simplex).

One then constructs a \emph{barcode diagram} from the persistence module, a visual representation from which one can deduce the homology of the underlying manifold, given a sufficiently good point sample. Construction of the barcode diagram is possible due to the structure theorem \cite{interval_decomp_source}, which says that a persistence module can be decomposed into a direct sum of simple components known as \emph{interval modules}.
The other important theorem that makes persistent homology interesting is the stability theorem, which guarantees the stability of the barcode diagram under small perturbations of the initial data.

Despite the stability theorem, persistent homology suffers from robustness issues. For example, persistent homology is highly unstable to outliers in the point data \cite{BotnanLesnickMultipersistence}.
A possible solution to this is \emph{multipersistent homology} \cite{multipersistence_source_1, multipersistence_source_2}, where instead of constructing an $\R_{\ge 0}$--parameterized filtration from the point data, one constructs an $(\R_{\ge 0})^k$--parameterized filtration for some integer $k$. This would allow one to, for example, take into account the certainty of point data and in this way avoid the outlier problem. Analogously to the single-parameter case, multipersistent homology gives a \emph{multipersistence module}, i.e., a functor
\begin{equation*}
    F \colon (\nnR)^k \to \Vecs_{\Fb}.
\end{equation*}
However, an analogue of the structure theorem fails to hold for such functors when $k \ge 2$. Thus, one cannot construct a barcode diagram for multipersistent homology, which is a well-known problem in TDA.

Instead of fixing the homology degree $i$ when constructing a (multi-)persistence module, we could instead consider the functor $H_{\bullet}(-; \Fb)$ that sends a simplicial complex $X$ to the chain complex with $H_i(X;\Fb)$ in degree $i$ and only zero differentials. This would give a functor
\begin{equation*}
    F \colon (\nnR)^k \to \Ch_{\Fb},
\end{equation*}
where $\Ch_{\Fb}$ is the category of chain complexes over $\Vecs_{\Fb}$.
As we care about such functors up to quasi-isomorphism, we can study them from the point of view of model categories, using the injective model structure on $\Ch_{\Fb}$ \cite[Theorem 2.3.13]{hovey}.

\subsubsection{Functor calculus}

Functor calculus is a general framework for studying functors, and is inspired by calculus of functions, and in particular the Taylor expansion of a function. The general idea of functor calculus is to introduce the notion of a \emph{degree $n$ functor}, analogous to a degree $n$ polynomial, and then study the \emph{degree $n$ approximations} of a functor. Typically, one gets a tower of approximations of the following form, which under reasonable conditions converges in a well-defined sense to the functor studied.
\begin{equation}\label{eq:calculus_tower}
\begin{aligned}
\begin{tikzpicture}
    \diagram{d}{3em}{3em}{
        \  & \vdots \\
        \ & P_1 F \\
        F & P_0 F \\
    };
    
    \path[->,font = \scriptsize, midway]
    (d-2-2) edge (d-3-2)
    (d-1-2) edge (d-2-2)
    (d-3-1) edge (d-1-2)
    (d-3-1) edge (d-2-2)
    (d-3-1) edge (d-3-2);
\end{tikzpicture}
\end{aligned}
\end{equation}

Among the flavors of functor calculus are Goodwillie calculus \cite{goodwillie_calculus, goodwillie_calculus2, goodwillie_calculus3}, abelian functor calculus \cite{abelian_calculus} and manifold calculus \cite{Goodwillie_1999}. Manifold calculus is particularly relevant to our project, as it is used to study a specific class of functors from a poset to spaces.

%The aim of manifold calculus is to study embeddings of one manifold into another. This is done by studying a specific contravariant functor on the poset of open subsets of a manifold. Explicitly, given manifolds $M$ and $N$, one studies the functor $F$ that sends an open subset $V$ of $M$ to the space of embeddings $\textrm{emb}(V, N)$. Manifold calculus then gives a tower of approximations to $F$.

Manifold calculus has been used successfully to study the homotopy type of components of the space of embeddings $\textrm{emb}(M, N)$ for manifolds $M, N$. For example, manifold calculus has been used to study long knots, which are embeddings from $\mathbb{R}$ into $\mathbb{R}^n$ that are linear outside a compact set. Using manifold calculus to study these knots is particularly effective when $n > 3$, as in this case the tower of approximations for $\textrm{emb}(M, N)$ converges \cite{Goodwillie_1999}. In \cite{volic2006calculus}, it is explained how, when $n > 3$, manifold calculus gives rise to a spectral sequence that can be used to determine the rational homotopy type of the space of long knots into $\mathbb{R}^n$.

There is also a dual notion of functor calculus, called \emph{functor cocalculus}. Functor cocalculus, produces a \emph{telescope} of approximations, i.e., a tower with all the arrows inverted. An example of functor cocalculus is McCarthy's dual calculus \cite{dual_calculus}.

\subsubsection{Related work}

Relations between multipersistence modules and functor calculus have been explored previously in \cite[Chapter 5.7]{lerch_thesis}. In \cite{lerch_thesis}, however, the author uses the existing theory of Goodwillie calculus and applies it to functors between $\infty$--categories.
Our work is, as far as we know, the first attempt at creating a flavor of functor calculus for functors from a poset, and to apply such a framework to study multipersistence modules.

We are also aware of another ongoing research project on this topic, by Nicolas Berkouk and Grégory Ginot, that is yet to be published. While that project also treats functor calculus and multipersistence modules, it has both a different goal and a fundamentally different approach than our article.

\subsection{Main contributions}

One of the most universal concepts in applied mathematics is the notion of approximation. Approximations allow us to replace complicated objects with objects that are easier to interpret or to compute, and both finding and understanding such approximations is a central theme in mathematics.
We introduce in this article \emph{poset cocalculus}, a flavor of functor cocalculus, as a tool for studying approximations in topological data analysis.
We further give several examples where poset cocalculus gives meaningful approximations when applied to objects that are relevant in TDA.
The versatility that these examples display leads us to believe that poset cocalculus can become a useful tool in TDA, both for finding new approximations and for helping us better understand the approximations that are already used.

In another preprint, \cite{hemFunctorCalculusMultipersistence}, we apply the theory developed here to study multipersistence modules. In particular, we show that codegree 1 bipersistence modules are interval decomposable, and establish a deep connection between the notion of \emph{middle-exact} multipersistence modules and poset cocalculus.

Poset cocalculus studies functors from a lattice (i.e., a poset where pairwise supremum and infimum are always defined) to a model category.
We define a functor $F \colon P \to \M$ to be \emph{codegree $n$} if it sends strongly bicartesian $(n+1)$--cubes to homotopy cocartesian $(n+1)$--cubes (see \autoref{sec:cubical_diags} for the definitions of these notions).
In \autoref{sec:taylor_telescope}, we define a notion of \emph{dimension} for elements in a lattice. We then define the codegree $n$ approximation of a functor $F$, denoted $T_n F$, as a homotopy left Kan extension of $F$ restricted to the elements of join-dimension $\le n$. This gives a telescope, called the \emph{Taylor telescope},
\begin{equation*}
\begin{aligned}
\begin{tikzpicture}
    \diagram{d}{3em}{3em}{
        \vdots & \ \\
        T_1 F & \ \\
        T_0 F & F. \\
    };
    
    \path[->,font = \scriptsize, midway]
    (d-2-1) edge (d-1-1)
    (d-1-1) edge (d-3-2)
    (d-3-1) edge (d-2-1)
    (d-2-1) edge (d-3-2)
    (d-3-1) edge (d-3-2);
\end{tikzpicture}
\end{aligned}
\end{equation*}

We then prove the following two theorems, which together say that this is a good notion of codegree $n$ approximation, under certain restrictions on the source lattice.
\begin{thmx}\label{theoremA_introduction}
    If $P$ is a distributive lattice, then for every functor $F \colon P \to \M$ to a good model category, $T_k F$ is codegree $k$.
\end{thmx}
\begin{thmx}\label{theoremB_introduction}
    Let $P$ be a join-factorization lattice, and let $F \colon P \to \M$ be a functor to a good model category. If $F$ is codegree $n$, then $T_n F \simeq F$.
\end{thmx}
The definitions of distributive lattices and join-factorization lattices are given in \autoref{sec:deg_n_defns} and \autoref{sec:taylor_telescope}.
Both results apply to a large class of lattices, including $(\nnR)^n$, $(\Nn, |)$ (i.e., the positive integers under division), and $(\Po(V), \subseteq)$ for $V$ finite.
We also show that the Taylor telescope of $F \colon P \to \M$ converges to $F$ whenever $P$ is a join-factorization lattice.
Finally, we prove that for any functor $F \colon P \to \M$ from a join-factorization lattice, the map $T_n F \to F$ is the universal map from a codegree $n$ functor (in the homotopy category).

% TODO: rewrite the following two paragraphs?
% We further give several examples that illustrate the utility of poset cocalculus, and in particular \autoref{theoremA_introduction} and \autoref{theoremB_introduction}.
%We compute the codegree 1 approximation of multipersistence modules, both when viewed as functors $F \colon (\nnR)^n \to \Vecs_{\Fb}$ (where $\Vecs_{\Fb}$ is equipped with the trivial model structure), and when viewed as functors $F \colon (\nnR)^n \to \Ch_{\Fb}$ (equipped with the injective model structure).
%We demonstrate that the codegree 1 approximation of a bipersistence module is interval decomposable.

Furthermore, we prove that multipersistence modules are stable under codegree $n$ approximation, given an appropriate notion of interleaving distance. This result can be summed up in the following inequalities.
\begin{equation*}
    d^{\Lambda}_H (T_1 F, T_1 G) \le d^{\Lambda}_H (T_2 F, T_2 G) \le \dots \le d^{\Lambda}_H (F, G)
\end{equation*}
This stability result opens up the possibility for using codegree $n$ approximations as a step in a TDA pipeline. This could be particularly useful in multipersistence with many parameters, where computational cost is a problem. The codegree $n$ approximations are indeed easier to compute, as they are determined by the values of the module on a small subset of the source poset.

We further give several examples of poset cocalculus applied to functors where both the source and target category are posets. We show that both the concept of $n$--skeletal simplicial complexes and $n$--coskeletal simplicial complexes can be phrased in terms of poset cocalculus. We then prove the following result about filtrations.
\begin{proposition*}
    Let $V \subset \R^d$ be a finite point cloud, and let $f_{\mathbf{C}} \colon \Po(V) \to \nnR$ and $f_{\mathbf{VR}} \colon \Po(V) \to \nnR$ be the Čech filtration function and the Vietoris-Rips filtration function, respectively.
    Then,
    \begin{equation*}
        f_{\mathbf{VR}} = 2 \cdot T_2 f_{\mathbf{C}}.
    \end{equation*}
\end{proposition*}
This proposition gives a new perspective on the role of the Vietoris-Rips filtration as a less computationally expensive approximation of the Čech filtration.

% We also illustrate some connections between poset cocalculus and discrete Morse theory \cite{Forman1998}, which is a field with important applications to TDA, particularly its computational aspects (see, e.g., \cite{Bauer_2021}). We show how being a critical $n$--simplex is somehow a ``failure to be degree $n$'', and we then formulate the Morse complex in terms of poset cocalculus.
% This formulation is interesting as it could potentially lead to generalizations of discrete Morse theory, which in turn could lead to new computational tools in TDA.

Finally, we show how poset cocalculus has a dual version, \emph{poset calculus}, that gives rise to a Taylor tower instead of a Taylor telescope. We then give an example of a computation in poset calculus that illustrate the utility of working with model categories. Given two simplicial complexes $X$ and $Y$, we study the functor $F$ that sends a subcomplex $Z \subseteq X$ to the space $|\Hhom(Z, Y)|$, where $\Hhom(-,-)$ is the internal Hom functor in the category of simplicial complexes and $| - | $ denotes geometric realization. It turns out that the \emph{degree 1} approximation of $F$ (dual of the codegree 1 approximation) in poset calculus is the functor that sends $Z \subseteq X$ to the space $\Hhom(|Z|, |Y|)$. This can be interpreted as saying that the degree 1 approximation of the space of simplicial maps between simplicial complexes is the space of continuous maps between their realizations.

\subsection{Organization}

In section 2, we review the needed preliminaries on homotopy colimits, and we prove some results on cubical homotopy colimits. Section 3 and section 4 are where we define poset cocalculus and prove our main theorems. In section 3, we give the definition of a codegree $n$ functor, and we prove results that classify strongly bicartesian cubes in distributive lattices. In section 4, we first define the join-dimension of an element in a lattice. We then define codegree $n$ approximations and the Taylor telescope, and we prove several results on these constructions, including \autoref{theoremA_introduction}, \autoref{theoremB_introduction} and a convergence statement. The remaining four sections are dedicated to examples of poset cocalculus. In section 5, we give some basic examples, as well as non-examples. Section 6 deals with multipersistence, and in particular stability. In section 7, we give several examples with poset-valued functors, and explore connections to the Čech and Vietoris-Rips filtrations. Finally, in section 8, we introduce the dual \emph{poset calculus}, and we give an example that relates to the homotopy type of Hom simplicial complexes.

\subsection{Notation}

Given two elements $x, y$ of a poset, we denote their \emph{least upper bound}, or \emph{join}, with $x \vee y$ and their \emph{greatest lower bound}, or \emph{meet}, with $x \wedge y$. A \emph{lattice} is a poset in which $x \vee y$ and $x \wedge y$ are defined for all pairs of elements $(x, y)$.

Let $R$ be a ring. We denote by $\Ch_{R}$ the category of unbounded chain complexes over $R$--modules. Unless otherwise specified, we endow $\Ch_R$ with the \emph{injective model structure}, in which the weak equivalences are the objectwise weak equivalences and the cofibrations are the objectwise monomorphisms. For the existence of this model structure, see \cite[Theorem 2.3.13]{hovey}.

Let $\M$ be a model category. We say that $\M$ is a \emph{good} model category if $\Fun(P, \M)$ admits the projective model structure for all posets $P$. Examples of good model categories include
\begin{itemize}
    \item Cofibrantly generated model categories \cite[Theorem 11.6.1]{hirschhorn}. Examples of cofibrantly generated model categories include:
    \begin{itemize}
        \item $\Ch_R$ equipped with the injective model structure \cite[Theorem 2.3.13]{hovey}.
        \item Simplicial sets, with the classical (Kan) model structure \cite[Example 11.1.6]{hirschhorn}.
        \item (Cofibrantly generated) topological spaces, with the classical (Quillen) model structure \cite[Example 11.1.8]{hirschhorn}\cite[Theorem 2.4.25]{hovey}.
    \end{itemize}
    \item Given a category $\C$, the trivial model structure on $\C$ where the weak equivalences are the isomorphisms in $\C$ and the cofibrations and fibrations are both all morphisms. In this case, the projective model structure on $\Fun(P, \C)$ is simply the trivial model structure.
\end{itemize}

Let $I$ be a small category and $\M$ a good model category. We define the \emph{homotopy colimit} of a functor $F \colon I \to \M$ as
\begin{equation*}
    \underset{I}{\hocolim} \ F \ = \ \underset{I}{\colim} \ QF,
\end{equation*}
where $QF$ is a cofibrant replacement of $F \in \Fun(I, \M)$, and we consider $\Fun(I, \M)$ to be equipped with the projective model structure.

\subsection*{Acknowledgments}

I would like to thank my supervisor, Kathryn Hess, whose ideas, advice and support contributed greatly to this work.

\section{Preliminaries on homotopy colimits}

\subsection{Kan extensions}
Given functors $\alpha \colon I \to J$ and $F \colon I \to \C$, let $\Lan_{\alpha} F$ denote the left Kan extension of $F$ along $\alpha$.
When $I$ and $J$ are posets, this is given explicitly by % maybe TODO: refer to Emily Riehl's book
\begin{equation*}
    (\Lan_{\alpha} F)(x) = \underset{y \in J, \alpha(y) \le x}{\colim} F(y)
\end{equation*}
where the colimit is taken over the subposet of $J$ as indicated.
Furthermore,
\begin{equation*}
    (\Lan_{\alpha} F)(x \le x') \colon \underset{y \in J, \alpha(y) \le x}{\colim} F(y) \to \underset{y \in J, \alpha(y) \le x'}{\colim} F(y)
\end{equation*}
is the morphism induced by the colimit cone $\{F(y) \to \colim_{y \in J, \alpha(y) \le x'} F(y)\}$.

\subsection{Homotopy colimits and homotopy left Kan extensions}
For $I$ a small category and $\M$ a model category, the \emph{homotopy colimit} is the functor
\begin{equation*}
    \underset{I}{\hocolim} \colon \Fun(I, \M) \to \M
\end{equation*}
defined as the total left derived functor of $\colim_I \colon \Fun(I, \M) \to \M$, where we consider the weak equivalences in $\Fun(I, \M)$ to be the objectwise weak equivalences.

Likewise, given a model category $\M$ and a functor of small categories $\alpha \colon I \to J$, the \emph{homotopy left Kan extension} is the functor
\begin{equation*}
    \hoLan_{\alpha} \colon \Fun(I, \M) \to \Fun(J, \M)
\end{equation*}
defined as the total left derived functor of the left Kan extension $\Lan_{\alpha} \colon \Fun(I, \M) \to \Fun(J, \M)$.

The existence of the homotopy colimits and homotopy left Kan extensions is guaranteed by part 2 and 3, respectively, in \cite[Theorem 11.3]{Scherer}. 

\subsection{The projective model structure}

Under certain conditions, such as when $\M$ is a cofibrantly generated model category, we can define the \emph{projective model structure} on $\Fun(I, \M)$. This is a model structure where
\begin{itemize}
    \item the weak equivalences are the natural transformations that are objectwise weak equivalences,
    \item the fibrations are the natural transformations that are objectwise fibrations,
    \item the cofibrations are the natural transformations with the left lifting property with respect to acyclic fibrations.
\end{itemize}

When $\Fun(I, \M)$ has the projective model structure, the adjunction
\[
    \begin{tikzcd}
        \Fun(I, \M) \arrow[r, shift left=1ex, "\colim_I"{name=G, yshift=1pt}] & \M \arrow[l, shift left=.5ex, "\const"{name=F}]
        \arrow[phantom, from=F, to=G, , "\scriptscriptstyle\boldsymbol{\bot}"],
    \end{tikzcd}
\]
is a Quillen adjunction. Hence, for a functor $F \colon I \to \M$, we can compute $\hocolim_I F$ by choosing a cofibrant replacement $QF$ of $F$ and computing $\colim_I (QF)$.

The same principle holds for homotopy left Kan extensions. When $\Fun(I, \M)$ and $\Fun(J, \M)$ both admit the projective model structure, then the adjunction
\[
    \begin{tikzcd}[column sep=large]
        \Fun(I, \M) \arrow[r, shift left=1ex, "\Lan_{\alpha}"{name=G, yshift=1pt}] & \Fun(J, \M) \arrow[l, shift left=.5ex, "\alpha^{*}"{name=F}]
        \arrow[phantom, from=F, to=G, , "\scriptscriptstyle\boldsymbol{\bot}"],
    \end{tikzcd}
\]
is a Quillen adjunction (because if $\eta$ is an objectwise (acyclic) fibration, then so is $\alpha^*(\eta)$). Hence, for a functor $F \colon I \to \M$ with cofibrant replacement $QF$,
\begin{equation*}
    \hoLan_{\alpha} F \simeq \Lan_{\alpha} (QF). 
\end{equation*}
Note also that, because $(\Lan_{\alpha} \dashv \alpha^{*})$ is a Quillen adjunction, when $F \colon I \to \M$ is a cofibrant functor, $\Lan_{\alpha} F$ is also a cofibrant functor.

\begin{lemma}\label{lemma:hocolim_hlke_comp}
    Let $\M$ be a model category, and let $\alpha \colon I \to J$ be a functor between small categories. Suppose that $\Fun(I, \M)$ and $\Fun(J, \M)$ both have the projective model structure.
    Then, for any $F \colon I \to \M$,
    \begin{equation*}
        \underset{J}{\hocolim} \left( \underset{\alpha}{\hoLan} \ F \right) \simeq \underset{I}{\hocolim} \ F.
    \end{equation*}
\end{lemma}
\begin{proof}
    The functor $\colim_J \circ \Lan_{\alpha}$ is the left adjoint of $\alpha^{*} \circ \const_J = \const_I$, and hence is naturally isomorphic to $\colim_I$. Now, let $F \colon I \to \M$, and let $QF$ be a cofibrant replacement of $F$. Then,
    \begin{align*}
        \underset{J}{\hocolim} \left( \underset{\alpha}{\hoLan} \ F \right) &= \underset{J}{\hocolim} \left( \underset{\alpha}{\Lan} \ QF \right) \\
        & \simeq \underset{J}{\colim} \left( \underset{\alpha}{\Lan} \ QF \right) \\
        & \cong \underset{I}{\colim} \ QF \\
        & = \underset{I}{\hocolim} \  F, \\
    \end{align*}
    where we use that $\Lan_{\alpha} QF$ is cofibrant because $QF$ is cofibrant and both functor categories are equipped with the projective model structure.
\end{proof}

\subsection{Homotopy final functors}

\begin{definition}\label{def:homotopy_finality}
    A functor $F \colon I \to J$ between small categories is \emph{homotopy final} if the under-category $x \downarrow F$ has contractible nerve for every $x \in J$.
\end{definition}

The following is Theorem 30.5 in \cite{Scherer}.
\begin{proposition}\label{prop:homotopy_finality}
    Let $I$ and $J$ be small categories, let $\M$ be a model category, and let $F \colon J \to \M$ be a functor. Suppose $\alpha \colon I \to J$ is a homotopy final functor. Then,
    \begin{equation*}
        \underset{I}{\hocolim} (F \circ \alpha) \simeq \underset{J}{\hocolim} F.
    \end{equation*}
\end{proposition}

\subsection{The Reedy model structure}

When our source category $I$ admits certain structure, we can define the \emph{Reedy model structure} on $\Fun(I, \M)$ for an arbitrary model category $\M$. Furthermore, when $I$ is a \emph{directed category}, we can define a Reedy model structure on $\Fun(I, \M)$ where the fibrations are objectwise fibrations and weak equivalences are objectwise weak equivalences. In other words, we can give $\Fun(I, \M)$ a Reedy model structure that coincides with the projective model structure.

Examples of directed categories include posets that satisfy the \emph{descending chain condition}.
\begin{definition}\label{def:descending_chain}
    A poset $P$ satisfies the \emph{descending chain condition} if every nonempty subset of $P$ has a minimal element.
\end{definition}
In particular, finite posets satisfy the descending chain condition.
Hence, when $P$ is a finite poset and $\M$ is any model category, we can equip $\Fun(P, \M)$ with the projective model structure. In this case, the weak equivalences are the objectwise weak equivalences, the fibrations are the objectwise fibrations, and the cofibrations are the natural transformations $F \to G$ such that for all $x \in P$, the morphism
\begin{equation}\label{eq:reedy_poset_cofibration}
    \underset{y < x}{\colim} \ G(y) \coprod_{\underset{y < x}{\colim} \ F(y)} F(x) \quad \to \quad G(x)
\end{equation}
is a cofibration.
In particular, the cofibrant functors are the functors $F$ such that
\begin{equation}\label{eq:reedy_poset_cofibrant_ob}
    \underset{y < x}{\colim} \ F(y) \to F(x)
\end{equation}
is a cofibration for all $x \in P$.

% Sources for this are: https://ncatlab.org/nlab/show/Reedy+model+structure , https://ncatlab.org/nlab/show/Reedy+category#direct_and_inverse_categories

It follows from \eqref{eq:reedy_poset_cofibration} that restricting a Reedy cofibration to a lower set (i.e., a subposet $L \subseteq P$ such that whenever $x \in L$ and $y \le x$ we have $y \in L$) yields a Reedy cofibration. In particular, restricting to a lower set preserves cofibrant functors in the Reedy model structure.

\subsection{Telescopes}

A \emph{telescope} in a model category $\M$ is a functor from the poset $\Nn$ to $\M$. As $\Nn$ satisfies the descending chain condition, we can equip $\Fun(\Nn, \M)$ with the projective model structure, and thus we can always define the homotopy colimit of a telescope.

\begin{lemma}\label{lemma:telescope_hocolim}
    Let $\M$ be a model category, let $T \colon \Nn \to \M$ be a telescope, and let $\zeta \colon T \to c$ be a cocone to an object in $\M$. Suppose that there exists an $N$ such that $\zeta_m \colon T(m) \to c$ is a weak equivalence for all $m \ge N$. Then, the map
    \begin{equation*}
        \hocolim_{\Nn} T \to c
    \end{equation*}
    is a weak equivalence.
\end{lemma}

\begin{proof}
First, observe that applying 2-out-of-3 on
\begin{equation*}
    \zeta_{N+i} = \zeta_{N+i} \circ T(N \le N+i),
\end{equation*}
we see that $T(N \le N+i)$ is a weak equivalence for all $i \ge 0$.

Let $q \colon QT \to T$ be a cofibrant replacement of $T$ so that $\hocolim_{\Nn} T = \colim_{\Nn} QT$.
Let $\hat QT$ be the telescope
\begin{center}
\begin{tikzcd}[column sep = 1.8em]
QT(0) \arrow[r] & QT(1) \arrow[r] & \cdots \arrow[r] & QT(N) \arrow[r, "="] & QT(N) \arrow[r, "="] & QT(N) \arrow[r] & \cdots,
\end{tikzcd}
\end{center}
and let $p \colon \hat QT \to QT$ be the natural transformation
\begin{center}
\begin{tikzcd}[column sep = 2.2em]
QT(0) \arrow[d, "="] \arrow[r] & \cdots \arrow[r] & QT(N) \arrow[r, "="] \arrow[d, "="] & QT(N) \arrow[r, "="] \arrow[d, "QT(N\le N+1)"] & QT(N) \arrow[r] \arrow[d, "QT(N\le N+2)"] & \cdots \\
QT(0) \arrow[r]                   & \cdots \arrow[r] & QT(N) \arrow[r] & QT(N+1) \arrow[r]               & QT(N+2) \arrow[r]                                   & \cdots.
\end{tikzcd}
\end{center}
The telescope $\hat Q T$ is cofibrant by \autoref{eq:reedy_poset_cofibrant_ob}. Furthermore, the vertical maps are weak equivalences, by 2-out-of-3 on
\begin{equation*}
    q_{N+i} \circ QT(N \le N+i) = T(N \le N+i) \circ q_N.
\end{equation*}
Thus, $p$ is a weak equivalence between cofibrant functors, so
\begin{equation*}
    \colim_{\Nn} QT \simeq \colim_{\Nn} \hat QT \simeq QT(N) \simeq c.
\end{equation*}
%Thus, $p$ is a weak equivalence between cofibrant functors, so the map
%\begin{equation*}
%    QT(N) = \colim_{\Nn} \hat QT \ \to \ \colim_{\Nn} QT
%\end{equation*}
%is a weak equivalence.
%
%Finally, applying 2-out-of-3 on
%\begin{center}
%\begin{tikzcd}
%QT(N) \arrow[rr] \arrow[rd] &   & \colim_{\Nn} QT \arrow[ld] \\
%                            & c &                        
%\end{tikzcd}
%\end{center}
%concludes the proof.
\end{proof}

\subsection{Cubical diagrams}
\label{sec:cubical_diags}

Let $[k] = \{0, \dots, k\}$, and let $\Po_{k+1}$ be the power set of $[k]$ (which is a poset ordered by inclusion). 
For $C$ any category, we call functors $\X \colon \Po_k \to C$ \emph{$k$--cubes} in $C$. 

\begin{definition}
    Let $C$ be a category.
    A $(k+1)$--cube $\X \colon \Po_{k+1} \to C$ is \emph{cocartesian} if the canonical map
    \begin{equation*}
        \underset{S \subsetneq [k]}{\colim \X(S)} \to \X([k])
    \end{equation*}
    is an isomorphism.

    Similarly, $\X$ is \emph{cartesian} if the canonical map
    \begin{equation*}
        \X(\emptyset) \to \underset{S \subseteq [k], S \neq \emptyset}{\lim \X(S)}
    \end{equation*}
    is an isomorphism.
\end{definition}

% possible example: pullbacks and pushouts

\begin{definition}
    Let $C$ be a category.
    A $(k+1)$--cube $\X \colon \Po_{k+1} \to C$ is \emph{strongly cartesian} (resp. \emph{strongly cocartesian}) if each face of dimension $\ge 2$ is cartesian (resp. cocartesian).

    If $\X$ is both strongly cartesian and strongly cocartesian, it is called \emph{strongly bicartesian}.
\end{definition}

\begin{remark}\label{rmk:2_faces_sufficient}
    For any cube, if all faces of dimension 2 are (co)cartesian, then so are all faces of dimension $> 2$. Hence, a cube is strongly (co)cartesian if all faces of dimension 2 are (co)cartesian. For more details, see  \cite[p.\ 272]{cubicalhtpy}. (This only covers the case where the target category is topological spaces, however the proof is straightforward and works in any category).
\end{remark}

\begin{definition}
    Let $\M$ be a model category.
    A $(k+1)$--cube $\X \colon \Po_{k+1} \to \M$ is \emph{homotopy cocartesian} if the canonical map
    \begin{equation}\label{eq:hocolim_XS}
        \underset{S \subsetneq [k]}{\hocolim \X(S)} \to \X([k])
    \end{equation}
    is a weak equivalence.

    Similarly, $\X$ is \emph{homotopy cartesian} if the canonical map
    \begin{equation*}
        \X(\emptyset) \to \underset{S \subseteq [k], S \neq \emptyset}{\holim \X(S)}
    \end{equation*}
    is a weak equivalence.
\end{definition}

%\begin{remark}
%When $\C$ is a model category with basepoint 0, the map in \eqref{eq:hocolim_XS} is a weak equivalence if and only if
%\begin{equation*}
%    \thocofib \X \simeq 0.
%\end{equation*}
%See \cite[Chapter 5.9]{cubicalhtpy} for more details.
%\end{remark}

\begin{lemma}\label{lemma:cube_hocolim_as_pushout}
    Let $\M$ be a model category. Given a $(k+1)$--cube $\X \colon \Po_{k+1} \to \M$, with $k \ge 1$, there is a weak equivalence
    \begin{equation*}
        \underset{S \subsetneq [k]}{\hocolim} \ \X(S) \ \simeq \ \hocolim \left( \X([k-1]) \leftarrow \underset{S \subsetneq [k-1]}{\hocolim} \X(S) \to \underset{S \subsetneq [k-1]}{\hocolim} \X(S \cup \{k\}) \right).
    \end{equation*}
\end{lemma}
\begin{proof}
As $\Po_{k+1}$ is a finite poset, $\Fun(\Po_{k+1}, \M)$ admits the projective model structure. Take a cofibrant replacement $Q\X$ of $\X$.
We want to show that
\begin{align*}
    &\hocolim \left( \X([k-1]) \leftarrow \underset{S \subsetneq [k-1]}{\hocolim} \X(S) \to \underset{S \subsetneq [k-1]}{\hocolim} \X(S \cup \{k\}) \right) \\
    \simeq \ & \colim \left( Q\X([k-1]) \leftarrow \underset{S \subsetneq [k-1]}{\colim} Q\X(S) \to \underset{S \subsetneq [k-1]}{\colim} Q\X(S \cup \{k\}) \right). \\
\end{align*}
This is true if
\begin{enumerate}
    \item[(i)] the restriction $Q\X |_{P_1}$ is Reedy cofibrant, where $P_1 = \{S : S \subsetneq [k-1]\}$,
    \item[(ii)] the morphism $Q\X([k-1]) \leftarrow \underset{S \subsetneq [k-1]}{\colim} Q\X(S)$ is a cofibration in $\M$, and
    \item[(iii)] the induced natural transformation $\eta \colon Q\X|_{P_1} \to Q\X( - \cup \{k\})|_{P_1}$ is a Reedy cofibration.
\end{enumerate}
Now, (i) follows from the fact that $P_1$ is a lower set, and (ii) follows directly from \eqref{eq:reedy_poset_cofibrant_ob}. To show that (iii) is true, we need to show that \eqref{eq:reedy_poset_cofibration} is satisfied, i.e., that for each $S \subsetneq [k-1]$ the morphism
\begin{equation*}
    \underset{T \subsetneq S}{\colim} \ Q\X(T \cup \{k\}) \coprod_{\underset{T \subsetneq S}{\colim} \ Q\X(T)} Q\X(S) \quad \to \quad Q\X(S \cup \{k\})
\end{equation*}
is a cofibration.
This is precisely the morphism
\begin{equation*}
    \underset{T \subsetneq (S \cup \{k\})}{\colim} \ Q\X(T) \quad \to \quad Q\X(S \cup \{k\})
\end{equation*}
which is a cofibration because $Q\X$ is Reedy cofibrant.

Finally,
\begin{align*}
    \simeq \ & \colim \left( Q\X([k-1]) \leftarrow \underset{S \subsetneq [k-1]}{\colim} Q\X(S) \to \underset{S \subsetneq [k-1]}{\colim} Q\X(S \cup \{k\}) \right) \\
    \simeq \ & \underset{S \subsetneq [k]}{\colim} \  Q\X(S) \\
    \simeq \ & \underset{S \subsetneq [k]}{\hocolim} \  \X(S).
\end{align*}
\end{proof}

A map of $k$--cubes $\X \to \Yc$ can be considered as a single $(k+1)$--cube $\Zc$, where
\begin{equation*}
    \Zc (S) = 
    \begin{cases}
        \X(S), &\quad k \notin S, \\
        \Yc(S \setminus \{k\}), &\quad k \in S. \\
    \end{cases}
\end{equation*}
We now prove some properties of such maps of cubes.

\begin{lemma}\label{lemma:map_of_hococart_cubes}
    Let $\X \to \Yc$ be a map of $k$--cubes, with $k \ge 1$, and consider the $(k+1)$--cube $\Zc = \X \to \Yc$. If both $\X$ and $\Yc$ are homotopy cocartesian, then so is $\Zc$.
\end{lemma}
\begin{proof}
By \autoref{lemma:cube_hocolim_as_pushout}, we can rewrite $\hocolim_{S \subsetneq [k]} \Zc(S)$ as the homotopy pushout of
\begin{equation*}
    \X([k-1]) \leftarrow \underset{S \subsetneq [k-1]}{\hocolim} \X(S) \to \underset{S \subsetneq [k-1]}{\hocolim} \Yc(S).
\end{equation*}
Denote this homotopy pushout by $P$ and consider the diagram
\begin{center}
    \begin{tikzpicture}
        \diagram{d}{5em}{3em}{
            \underset{S \subsetneq [k-1]}{\hocolim} \X(S) & \underset{S \subsetneq [k-1]}{\hocolim} \Yc(S) & \ \\
            \X([k-1]) & P & \Yc([k-1]) \\
        };
        
        \path[->,font = \scriptsize, midway]
        (d-1-1)+(0,-15pt) edge node[midway, right]{$q$} (d-2-1)
        (d-1-1) edge node[midway, above]{$g$} (d-1-2)
        (d-2-1) edge node[midway, above]{$f$} (d-2-2)
        (d-1-2)+(0,-15pt) edge node[midway, right]{$p$} (d-2-2)
        (d-2-2) edge node[midway, above]{$\gamma$} (d-2-3);
    \end{tikzpicture}
\end{center}
% These things are maybe easier to see if I rather work with the colimit of a cofibrant replacement diagram.
As $\X$ is homotopy cocartesian, $q$ is a weak equivalence. Hence, $p$ is a weak equivalence. Furthermore, as $\Yc$ is homotopy cocartesian, $\gamma \circ p$ is a weak equivalence. By 2-out-of-3, $\gamma$ is a weak equivalence, as desired.
\end{proof}

\begin{lemma}\label{lemma:map_of_cubes_we}
    Let $\X \to \Yc$ be a map of $k$--cubes, with $k \ge 1$, and consider the $(k+1)$--cube $\Zc = \X \to \Yc$. If the map $\X \to \Yc$ is an objectwise weak equivalence, then $\Zc$ is homotopy cocartesian.
\end{lemma}
\begin{proof}
We again use \autoref{lemma:cube_hocolim_as_pushout} and rewrite $\hocolim_{S \subsetneq [k]} \Zc(S)$ as the homotopy pushout of
\begin{equation*}
    \X([k-1]) \leftarrow \underset{S \subsetneq [k-1]}{\hocolim} \X(S) \to \underset{S \subsetneq [k-1]}{\hocolim} \Yc(S).
\end{equation*}
Denote this homotopy pushout by $P$ and consider the diagram
\begin{center}
    \begin{tikzpicture}
        \diagram{d}{5em}{3em}{
            \underset{S \subsetneq [k-1]}{\hocolim} \X(S) & \underset{S \subsetneq [k-1]}{\hocolim} \Yc(S) & \ \\
            \X([k-1]) & P & \Yc([k-1]) \\
        };
        
        \path[->,font = \scriptsize, midway]
        (d-1-1)+(0,-15pt) edge node[midway, right]{$q$} (d-2-1)
        (d-1-1) edge node[midway, above]{$g$} (d-1-2)
        (d-2-1) edge node[midway, above]{$f$} (d-2-2)
        (d-1-2)+(0,-15pt) edge node[midway, right]{$p$} (d-2-2)
        (d-2-2) edge node[midway, above]{$\gamma$} (d-2-3);
    \end{tikzpicture}
\end{center}
As $\X \to \Yc$ is an objectwise weak equivalence (i.e., a Reedy weak equivalence), $g$ is a weak equivalence. Hence, $f$ is a weak equivalence. Furthermore, $\gamma \circ f \colon \X([k-1]) \to \Yc([k-1])$ is a weak equivalence. By 2-out-of-3, $\gamma$ is a weak equivalence, as desired.
\end{proof}

\section{Codegree $n$ functors on lattices}

\label{sec:deg_n_defns}

%From here on, all posets are assumed to have a minimal element, denoted 0. Hence, when we write \emph{lattice}, we mean \emph{lattice with a minimal element}.

\subsection{Codegree $n$ functors}

\begin{definition}\label{def:deg_n_bicartesian}
    Let $P$ be a lattice and $\M$ a good model category. A functor $F \colon P \to \M$ is \emph{codegree $n$} if it sends strongly bicartesian $(n+1)$--cubes to homotopy cocartesian cubes.

    In other words, $F$ is codegree $n$ if for each strongly bicartesian $(n+1)$--cube $\X \colon \Po_{n+1} \to P$, the canonical map
    \begin{equation}\label{eq:deg_n_bicartesian}
        \underset{S \subsetneq [n]}{\hocolim \left(F \circ \X(S) \right)} \to F \circ \X([n])
    \end{equation}
    is a weak equivalence.
\end{definition}

\begin{proposition}\label{prop:deg_increasing}
    Let $P$ be a lattice and $F \colon P \to \M$ a functor to a good model category. If $F$ is codegree $n$, then $F$ is codegree $m$ for all $m \ge n$.
\end{proposition}

\begin{proof}
    We proceed by induction. Suppose that $F$ is codegree $m$, and let $\Zc \colon \Po_{m+2} \to P$ be a strongly bicartesian $(m+2)$--cube.
    We can view the $\Zc$ as a map of $(m+1)$--cubes $\X \to \Yc$.
    Then $F \circ \X$ and $F \circ \Yc$ are homotopy cocartesian by the induction hypothesis.
    Thus, $F \circ \Zc = (F \circ \X \to F \circ \Yc)$ is also homotopy cocartesian, by \autoref{lemma:map_of_hococart_cubes}.
\end{proof}

\subsection{Pairwise covers}

\begin{definition}
    Let $P$ be a lattice, let $v \in P$ and let $k$ be a positive integer.
    A \emph{pairwise cover} of $v$ of size $k$ is a collection of elements $x^0, \dots, x^{k-1} \in P$, with $x^i \le v \ \forall i$, such that $x^i \vee x^j = v$ for all $i \neq j$.
\end{definition}

\begin{remark}
When $k \ge 2$, the fact that $x^i \vee x^j = v$ for $i \neq j$ implies $x^i \le v$ for all $i$. We add the requirement that $x^i \le v$ for all $i$ so that this also holds when $k = 1$ (in which case a pairwise cover is just any element $\le v$).
\end{remark}

\begin{definition}\label{def:poset_cube}
    Let $P$ be a lattice, and let $v \in P$. Let further $x^0, \dots, x^k$ be a pairwise cover of $v$.
    We define the $(k+1)$--cube 
    \begin{equation*}
        \X_{x^0, \dots, x^k} \colon \Po_{k+1} \to P
    \end{equation*}
    as follows.
    \begin{equation*}
        \X_{x^0, \dots, x^k} (S) = 
        \begin{cases}
            v, &\quad S = [k], \\ 
            \bigwedge_{i \notin S} x^i, &\quad \text{ otherwise.} \\
        \end{cases}
    \end{equation*}
\end{definition}

A lattice $P$ is \emph{distributive} if for all elements $x,y,z \in P$, we have the equality
\begin{equation}\label{eq:dist_lattice}
    x \wedge (y \vee z) = (x \wedge y) \vee (x \wedge z).
\end{equation}
Equivalently, a lattice $P$ is distributive if for all elements $x,y,z \in P$, we have
\begin{equation}\label{eq:dist_lattice_alt}
    x \vee (y \wedge z) = (x \vee y) \wedge (x \vee z)
\end{equation}
\cite[Chapter IX, Theorem 1]{Birkhoff}.

The following two lemmas show a correspondence between pairwise covers of size $k+1$ and strongly bicartesian $(k+1)$--cubes, in the case where the lattice is distributive.

\begin{lemma}\label{lem:bicartesian_from_codecomp}
    Let $P$ be a distributive lattice. For every $v \in P$ and every pairwise cover $x^0, \dots, x^k$ of $v$, the cube $\X_{x^0, \dots, x^k}$ is strongly bicartesian.
\end{lemma}

\begin{proof}
    As remarked in \autoref{rmk:2_faces_sufficient}, it suffices to show that each 2--dimensional face of the cube is bicartesian (i.e., both cartesian and cocartesian). These faces are either of the form
    \begin{center}
    \begin{tikzpicture}
        \diagram{d}{3em}{3em}{
            x^i \wedge x^j & x^i \\
            x^j & v \\
        };
        
        \path[->,font = \scriptsize, midway]
        (d-1-1) edge (d-2-1)
        (d-1-1) edge (d-1-2)
        (d-2-1) edge (d-2-2)
        (d-1-2) edge (d-2-2);
    \end{tikzpicture}
    \end{center}
    where $i \neq j$, or
    \begin{center}
    \begin{tikzpicture}
        \diagram{d}{3em}{3em}{
            \bigwedge_{l \notin S} x^l & \bigwedge_{l \notin S \cup \{i\}} x^l \\
            \bigwedge_{l \notin S \cup \{j\}} x^l & \bigwedge_{l \notin S \cup \{i, j\}} x^l \\
        };
        
        \path[->,font = \scriptsize, midway]
        (d-1-1) edge (d-2-1)
        (d-1-1) edge (d-1-2)
        (d-2-1) edge (d-2-2)
        (d-1-2) edge (d-2-2);
    \end{tikzpicture}
    \end{center}
    where $S \subset [k]$, and $i, j \in [k] \setminus S$ such that $i \neq j$.
    
    Both squares are clearly cartesian (because the upper left corner is the meet of the lower left and the upper right corner). The first square is cocartesian because $x^i \vee x^j = v$ from the definition of pairwise covers. The second square is cocartesian because
    \begin{align*}
        \left(\bigwedge_{l \notin S \cup \{i\}} x^l \right) \vee \left(\bigwedge_{l \notin S \cup \{j\}} x^l \right) &= \left(x^j \wedge \bigwedge_{l \notin S \cup \{i, j\}} x^l \right) \vee \left(x^i \wedge \bigwedge_{l \notin S \cup \{i, j\}} x^l \right) \\
        &= (x^j \vee x^i) \wedge \left( \bigwedge_{l \notin S \cup \{i, j\}} x^l \right) \\
        &= v \wedge \left( \bigwedge_{l \notin S \cup \{i, j\}} x^l \right) \\
        &= \bigwedge_{l \notin S \cup \{i, j\}} x^l, \\
    \end{align*}
    where the second equality follows from \eqref{eq:dist_lattice}, and the third equality follows from the fact that $i \neq j \Rightarrow x^i \vee x^j = v$.
\end{proof}

\begin{lemma}\label{lem:codecomp_from_bicartesian}
    % State that each strongly cocartesian cube gives a pairwise cover
    Let $P$ be a lattice, and let $\X \colon \Po_{k+1} \to P$ be a strongly bicartesian cube. For $i \in [k]$,\footnote{Recall that $[k] = \{0, \dots, k\}.$} let $x^i = \X\big([k] \setminus \{i\}\big)$. Then $x^0, \dots, x^k$ is a pairwise cover of $\X\big([k]\big)$.

    Furthermore, $\X = \X_{x^0, \dots, x^k}$, as defined in \autoref{def:poset_cube}.
\end{lemma}

\begin{proof}
    The fact that $x^0, \dots, x^k$ is a pairwise cover of $\X([k])$ follows from the fact that, for $i \neq j$, the following 2--face of $\X$ is cocartesian
    \begin{center}
    \begin{tikzpicture}
        \diagram{d}{3em}{3em}{
            \X\big([k] \setminus \{i, j\}\big) & \X\big([k] \setminus \{i\}\big) \\
            \X\big([k] \setminus \{j\}\big) & \X\big([k]\big), \\
        };
        
        \path[->,font = \scriptsize, midway]
        (d-1-1) edge (d-2-1)
        (d-1-1) edge (d-1-2)
        (d-2-1) edge (d-2-2)
        (d-1-2) edge (d-2-2);
    \end{tikzpicture}
    \end{center}
    meaning that $\X\big([k]\big) = \X\big([k] \setminus \{i\}\big) \vee \X\big([k] \setminus \{j\}\big) = x^i \vee x^j$.

    To show that $\X = \X_{x^0, \dots, x^k}$, we proceed by induction. Suppose that $\X\big([k] \setminus S\big) = \X_{x^0, \dots, x^n}\big([k] \setminus S\big)$ for sets $S \subseteq [k]$ such that $|S| \le n$. Let $T = \{t_1, \dots, t_{n+1}\} \subseteq [k]$. As the following 2--face is cartesian,
    \begin{center}
    \begin{tikzpicture}
        \diagram{d}{3em}{3em}{
            \X\big([k] \setminus \{t_1, \dots, t_{n+1}\}\big) & \X\big([k] \setminus \{t_1, \dots, t_n\}\big) \\
            \X\big([k] \setminus \{t_1, \dots, t_{n-1}, t_{n+1}\}\big) & \X\big([k] \setminus \{t_1, \dots, t_{n-1}\}\big), \\
        };
        
        \path[->,font = \scriptsize, midway]
        (d-1-1) edge (d-2-1)
        (d-1-1) edge (d-1-2)
        (d-2-1) edge (d-2-2)
        (d-1-2) edge (d-2-2);
    \end{tikzpicture}
    \end{center}
    it follows that
    \begin{align*}
        \X\big([k] \setminus & \{t_1, \dots, t_{n+1}\}\big) \\
        &= \X\big([k] \setminus \{t_1, \dots, t_{n-1}, t_{n+1}\}\big) \wedge \X\big([k] \setminus \{t_1, \dots, t_n\}\big)\\
        &= \X_{x^0, \dots, x^k}\big([k] \setminus \{t_1, \dots, t_{n-1}, t_{n+1}\}\big) \wedge \X_{x^0, \dots, x^k}\big([k] \setminus \{t_1, \dots, t_n\}\big)\\
        &= \left(\bigwedge_{i \notin ([k] \setminus \{t_1, \dots, t_{n-1}, t_{n+1}\})} x^i \right) \wedge \left(\bigwedge_{i \notin ([k] \setminus \{t_1, \dots, t_n\})} x^i \right) \\
        &= \left(\bigwedge_{i \notin ([k] \setminus \{t_1, \dots, t_{n+1}\})} x^i \right) \\
        &= \X_{x^0, \dots, x^k} \big([k] \setminus \{t_1, \dots, t_{n+1}\}\big).
    \end{align*}
    Hence, $\X(S) = \X_{x^0, \dots, x^k}(S)$ for all $S \subseteq [k]$.
    
\end{proof}

\begin{proposition}\label{prop:defns_consistent}
    Let $P$ be a distributive lattice and $\M$ a good model category. A functor $F\colon P \to \M$ is \emph{codegree $n$} if and only if the following condition holds. For every $v \in P$ and every size $(n+1)$ pairwise cover $x^0, \dots, x^n$ of $v$, the $(n+1)$--cube
    \begin{equation*}
        F \circ \X_{x^0, \dots, x^n} \colon \Po_{n+1} \to \M
    \end{equation*}
    is homotopy cocartesian.
\end{proposition}

\begin{proof}
    First, suppose $F$ is a codegree $n$ functor. By \autoref{lem:bicartesian_from_codecomp} the cube $\X_{x^0, \dots, x^n}$ is strongly bicartesian and hence $F \circ \X_{x^0, \dots, x^n}$ is homotopy cocartesian by \autoref{def:deg_n_bicartesian}.

    Now, suppose that $F \circ \X_{x^0, \dots, x^n}$ is homotopy cocartesian for any size $(n+1)$ pairwise cover $x^0, \dots, x^n$. Then, by \autoref{lem:codecomp_from_bicartesian}, $F \circ \X$ is homotopy cocartesian for any strongly bicartesian cube $\X$. Hence, $F$ is codegree $n$.
\end{proof}

\subsection{Examples}
From now on, whenever we write $P \times Q$ for posets $P$ and $Q$, we will assume that $P \times Q$ is equipped with the product order (i.e., $(x,y) \le (x',y')$ if and only if $x \le x'$ and $y \le y'$). For an element $v \in P_1 \times \dots \times P_n$, we will use $v_i$ to denote the $i$th coordinate of $v$, i.e., $v = (v_1, \dots, v_n)$.

Some examples of distributive lattices are the following.
\begin{itemize}
    \item Total orders. In this poset, join and meet are maximum and minimum, respectively.
    \item Arbitrary products of total orders. In this poset, join and meet are elementwise maximum and elementwise minimum, respectively.
    \item For any set $X$, the poset $\Po(X)$ where the partial order is given by inclusion. In this poset, join and meet are union and intersection, respectively. This poset is isomorphic to $\{0, 1\}^{\times |X|}$, i.e., a (possibly infinite) product of the 2-element total order. 
\end{itemize}

\begin{example}
A 1--cube consists of a single morphism between two objects $\X(0) \to \X(1)$. Any 1--cube is by definition strongly bicartesian. Hence, a functor $F \colon P \to \M$ (from a lattice to a good model category) is codegree 0 if and only if $F(x \le y)$ is a weak equivalence for every $x \le y$ in $P$.
\end{example}

\begin{remark}\label{rmk:one_x_is_v}
Let $P$ be a lattice, $\M$ a good model category and $F \colon P \to \M$ a functor.
Let $x^0, \dots, x^n$ be a pairwise cover of $v \in P$.
If $x^i = v$ for some $i$, then $F \circ \X_{x^0, \dots, x^n}$ is homotopy cocartesian, no matter what the functor $F$ is.
To see this, we make use of \autoref{lemma:map_of_cubes_we}.

Assume without loss of generality that $x^n = v$. We can view $\X_{x^0, \dots, x^n}$ as a map $\Yc \to \Zc$, where $\Yc$ and $\Zc$ are the $n$--cubes
\begin{equation*}
    \Yc (S) = 
    \begin{cases}
        x^n \wedge v, &\quad S = [n-1], \\ 
        x^n \wedge \bigwedge_{i \notin S} x^i, &\quad \text{otherwise.} \\
    \end{cases}
\end{equation*}
\begin{equation*}
    \Zc (S) = 
    \begin{cases}
        v, &\quad S = [n-1], \\ 
        \bigwedge_{i \notin S} x^i, &\quad \text{otherwise.} \\
    \end{cases}
\end{equation*}
Now, as $x^n = v$, we have that $x^n \wedge v = v$ and $x^n \wedge \bigwedge_{i \notin S} x^i$ = $\bigwedge_{i \notin S} x^i $, and hence the map $\Yc(S) \to \Zc(S)$ is the identity for all $S$.
Thus, $F \circ \Yc(S) \to F \circ \Zc(S)$ is also the identity for all $S$.
Hence, by \autoref{lemma:map_of_cubes_we}, $F \circ \X_{x^0, \dots, x^n} = F \circ \Yc \to F \circ \Zc$ is homotopy cocartesian.

We will repeatedly use this fact.
\end{remark}

\begin{example}\label{ex:deg_1_functor}
Consider the following functor from $\{0,1\}^2$ to $\Ch_{\Fb}$, where ${\Fb}$ in the diagram denotes the chain complex with the field ${\Fb}$ in degree 0 and 0 in all other degrees.
\begin{center}
\begin{tikzpicture}
    \diagram{d}{3em}{3em}{
        0 & \Fb \\
        0 & \Fb \\
    };
    
    \path[->,font = \scriptsize, midway]
    (d-1-1) edge (d-1-2)
    (d-2-2) edge node[right]{$1$} (d-1-2)
    (d-2-1) edge (d-1-1)
    (d-2-1) edge (d-2-2);
\end{tikzpicture}
\end{center}
This functor is codegree 1. We verify this by computing \eqref{eq:deg_n_bicartesian} of the cube $\X_{x^0, x^1}$ for various choices of $v \in P$ and pairwise covers $x^0, x^1$ of $v$. By \autoref{rmk:one_x_is_v}, we need only consider the choices where neither $x^0$ nor $x^1$ equals $v$, and this is only possible when $v = (1,1)$ and $\{x^0, x^1\} = \{(0,1), (1,0)\}$. It's easy to verify that the map in \eqref{eq:deg_n_bicartesian} is indeed a weak equivalence in this case.
\end{example}

\begin{example}\label{ex:deg_2_functor}
Consider the following functor from $\{0,1\}^2$ to $\Ch_{\Fb}$.
\begin{center}
\begin{tikzpicture}
    \diagram{d}{3em}{3em}{
        \Fb & \Fb \\
        0 & \Fb \\
    };
    
    \path[->,font = \scriptsize, midway]
    (d-1-1) edge node[above]{$1$} (d-1-2)
    (d-2-2) edge node[right]{$1$} (d-1-2)
    (d-2-1) edge (d-1-1)
    (d-2-1) edge (d-2-2);
\end{tikzpicture}
\end{center}
This functor is codegree 2, and it is not codegree 1. It is not codegree 1 because the cube defined by the pairwise cover $x^0=(0,1), x^1=(1,0)$ is not homotopy cocartesian (for that, we would need $F(1,1) \simeq \Fb \oplus \Fb$).

It is codegree 2 (again, using \autoref{rmk:one_x_is_v}) because any choice of $v \in P$ and size 3 pairwise cover $x^0, x^1, x^2$ of $v$ must satisfy $x^i = v$ for some $i$.
\end{example}

\begin{example}\label{ex:direct_sum_functor}
This functor is also codegree 1.
\begin{center}
\begin{tikzpicture}
    \diagram{d}{3em}{3em}{
        \Fb & \Fb \oplus \Fb \\
        0 & \Fb \\
    };
    
    \path[->,font = \scriptsize, midway]
    (d-1-1) edge node[above]{$\begin{pmatrix}1 \\ 0\end{pmatrix}$} (d-1-2)
    (d-2-2) edge node[right]{$\begin{pmatrix}0 \\ 1\end{pmatrix}$} (d-1-2)
    (d-2-1) edge (d-1-1)
    (d-2-1) edge (d-2-2);
\end{tikzpicture}
\end{center}
The computation is similar to that of the preceding examples.
\end{example}

\begin{example}\label{ex:Nn_2}
Consider any functor $F\colon \Nn \times \Nn \to \M$, where $\Nn \times \Nn$ is considered as a poset with the product order, and $\M$ is a good model category. Then $F$ is codegree $2$. To see this, let $v \in \Nn \times \Nn$, and let $x^0, x^1$ and $x^2$ be a pairwise cover of $v$. Write $x^i = (a_i, b_i)$ and $v = (u, w)$. It follows that 
\begin{align*}
    (u, w) &= \left(\max\{a_0, a_1\}, \max\{b_0, b_1\}\right)\\
    &= \left(\max\{a_0, a_2\}, \max\{b_0, b_2\}\right)\\
    &= \left(\max\{a_1, a_2\}, \max\{b_1, b_2\}\right)
\end{align*}
Hence, two of $\{a_0, a_1, a_2\}$ equal $u$ and two of $\{b_0, b_1, b_2\}$ equal $w$.
% TODO: This can be visualized
Thus, either $x_0$, $x_1$ or $x_2$ equals $v$, and $F \circ \X_{x^0, x^1, x^2}$ is homotopy cocartesian by \autoref{rmk:one_x_is_v}.
\end{example}

The previous example generalizes to $\Nn^{\times n}$ for any $n$. Furthermore, you can replace $\Nn$ with any totally ordered poset.

\begin{proposition}\label{thm:total_order_product}
If $P = P_1 \times \dots \times P_n$, where each $P_i$ is totally ordered, then any functor $F \colon P \to \M$ is codegree $n$.
\end{proposition}
\begin{proof}
The theorem follows from an argument similar to that in \autoref{ex:Nn_2}.

If $x \vee y = v$, then for each $i \in \{1, \dots, n\}$, either $x_i = v_i$ or $y_i = v_i$.
It follows that given any size $(n+1)$ pairwise cover $x^0, \dots, x^n$ of $v \in P$, then for all $i$, we have that $(x^j)_i = v_i$ for all but one $j$.
Hence, by the pigeonhole principle, one of the elements in the pairwise cover must equal $v$. The result follows from \autoref{rmk:one_x_is_v}.
\end{proof}

\begin{remark}
If $P = \R^n$, the proof of \autoref{thm:total_order_product} can be phrased in terms of linear algebra. If $x \vee y = v$, then for all $i \in \{1, \dots, n\}$, either $x_i = v_i$ or $y_i = v_i$. In other words, for all $i \in \{1, \dots, n\}$, either $v_i - x_i = 0$ or $v_i - y_i = 0$, which implies that $(v - x) \cdot (v - y) = 0$.

Thus, if you have $k$ elements $x^0, \dots, x^{k-1}$ such that $x^i \vee x^j$ for all $i \neq j$, then the vectors $(v - x^0), \dots, (v - x^{k-1})$ are pairwise orthogonal. Hence, if $k \geq n+1$, then one of these vectors must be 0.
\end{remark}

Observe that the functor in \autoref{ex:deg_1_functor} appears to be ``independent of one variable'', and thus behaves like a functor in only one variable. We formalize this into a general result.

\begin{proposition}\label{thm:independent_variables}
Let $P = P_1 \times \dots \times P_n$, where each $P_i$ is totally ordered, and let $F \colon P \to \M$ be a functor. Suppose further there is a set $I \subset \{1, \dots, n\}$ such that for any $i \in I$, $F(\id_{x_1}, \dots, \id_{x_{i-1}}, (x_i \le y_i), \id_{x_{i+1}}, \dots, \id_{x_n})$ is an isomorphism for all choices of $(x_1, \dots, x_n) \in P$ and $y_i \in P_i$.

Then $F$ has codegree $n - |I|$.
\end{proposition}
\begin{proof}
Let $P$ and $F$ be as in the proposition statement.
If two points $x = (x_1, \dots, x_n)$ and $v = (v_1, \dots, v_n)$, with $x \le v$, have the same values in coordinate $i$ for all $i \notin I$, then $F(x \le v)$ is an isomorphism. Furthermore, for any $a \in P$, $a \wedge x$ and $a \wedge v$ will also have the same values in coordinate $i$ for all $i \notin I$, and hence $F(a \wedge x \le a \wedge v)$ is an isomorphism.

Now, suppose we have $v \in P$ and $x^0, \dots, x^{n - |I|} \in P$ such that $x^i \wedge x^j = v$ for all $i \neq j$.
Then, by a similar counting argument as in \autoref{thm:total_order_product}, there must be an $x^i$ such that for all $j \notin I$, $(x^i)_j = v_j$.
The proposition now follows from a similar argument to that in \autoref{rmk:one_x_is_v}.
\end{proof}

% \begin{proposition}
% If $F$ and $G$ are cofibrant and codegree $n$, then the coproduct $F \coprod G$ has codegree $n$.
% \end{proposition}
% 
% \begin{proof}
% This follows from the fact that homotopy colimits commute, and that $F \coprod G$ is the homotopy coproduct of $F$ and $G$ when $F$ and $G$ are cofibrant.
% \end{proof}

\autoref{ex:direct_sum_functor} provides a functor that is a direct sum of two cofibrant codegree 1 functors.

\subsection{Relation to other work}

\subsubsection{Relation to manifold calculus}
The central definitions in poset cocalculus are inspired by the manifold calculus of Goodwillie and Weiss \cite{Weiss_1999, Goodwillie_1999}, where one also studies functors from posets arising naturally from manifolds. Poset cocalculus is not a strict generalization of the Goodwillie-Weiss theory, however. First, the definitions in poset cocalculus are dual: contravariant functors are replaced by covariant functors and homotopy cartesian cubes are replaced by homotopy cocartesian cubes. Second, in \cite{Weiss_1999}, Weiss restricts his attention to a subclass of functors that he calls \emph{good cofunctors}, while we make no such restriction, which leads to a different notion of (co)degree $n$ approximations.

We describe here the relation between degree $n$ functors in manifold calculus and codegree $n$ functors in poset cocalculus. We show that the definition of codegree $n$ in poset cocalculus can be viewed as a dual version of a generalization of the notion of degree $n$ in manifold calculus.

Given a manifold $M$, let $\bigo (M)$ be the poset consisting of open subsets of $M$, where we define a partial order relation by $U \le V$ if $U$ is a subset of $V$.
In manifold calculus, one studies contravariant functors
\begin{equation*}
    \bigo(M) \to \Spaces.
\end{equation*}

Given a contravariant functor $F\colon \bigo(M) \to \Spaces$, a set $V \in \bigo(M)$ and a collection of pairwise disjoint closed subsets $A_0, \dots, A_k$ of $V$, there is a $k+1$--cube $\X \colon \Po_{k+1} \to \Spaces$ given by
\begin{equation}\label{eq:cube}
    \X (S) = F\left(V \setminus \bigcup_{i \in S} A_i \right).
\end{equation}

\begin{definition}\label{def:degree_n_mf_calc}
A contravariant functor $F\colon \bigo (M) \to \Spaces$ is said to be of \emph{degree $n$} if for all $V$ and for all pairwise disjoint closed subsets $A_0, \dots, A_n$ of $V$, the cube $\X$ in \eqref{eq:cube} is homotopy cartesian.
\end{definition}
It follows that if $F$ is degree $n$, then it is degree $k$ for all $k \ge n$.

%We now an analogue of this concept for functors from an arbitrary poset $P$ (with suitable restrictions) to a model category.
We will compare this definition of degree $n$ functors to \autoref{def:deg_n_bicartesian}.
For two sets $U, V \in \bigo(M)$, their union $U \cup V$ equals the least upper bound of $U$ and $V$ in $\bigo(M)$ considered as a poset. This, in turn, equals the coproduct of $U$ and $V$ in $\bigo(M)$ considered as category. Similarly, the intersection $U \cap V$ corresponds to the greatest lower bound and to the categorical product. The complement operation, however, has no analogue in a general poset. We therefore reformulate the definitions above.

Let $F\colon M \to \Spaces$, $V \in \bigo(M)$, and $A_0, \dots, A_k$ be as above. Define $B_i = V \setminus A_i$. The $B_i$'s are open sets, and thus elements of $\bigo(M)$. The condition that the $A_i$'s are pairwise disjoint is equivalent to saying that for $i \neq j$, we have $B_i \cup B_j = V$.
Furthermore, $V \setminus \bigcup_{i\in S} A_i = \bigcap_{i \in S} (V \setminus A_i) = \bigcap_{i \in S} B_i$.
Hence, we can equivalently formulate \autoref{def:degree_n_mf_calc} as follows.
\begin{definition}\label{def:new_degree_n_mf_calc}
Let $F\colon \bigo (M) \to \Spaces$ be a contravariant functor. Then $F$ is said to be of \emph{degree $n$} if for all $V$ and for all open subsets $B_0, \dots, B_n$ of $V$, satisfying $B_i \cup B_j = V$ for $i \neq j$, the cube $\X$ defined by
\begin{equation*}
    \X (S) = 
    \begin{cases}
        F(V), &\quad S = \emptyset, \\ 
        F\left(\bigcap_{i \in S} B_i \right), &\quad \text{otherwise.} \\
    \end{cases}
\end{equation*}
is homotopy cartesian.
\end{definition}

We can now compare this to \autoref{def:deg_n_bicartesian}.
Let $P$ be a distributive lattice, $\M$ a good model category, and $F \colon P \to \M$ a (covariant) functor.
By \autoref{prop:defns_consistent}, $F$ is codegree $n$ (according to \autoref{def:deg_n_bicartesian}) if and only if the following condition holds.
\begin{itemize}
    \item 
    For every $v \in P$ and all $x^0, \dots, x^n$ such that $x^i \vee x^j = v$ for all $i \neq j$, the $(n+1)$--cube $\X$ defined by
    \begin{equation*}
        \X (S) = 
        \begin{cases}
            F(v), &\quad S = [n], \\ 
            F\left(\bigwedge_{i \notin S} x^i \right), &\quad \text{ otherwise.} \\
        \end{cases}
    \end{equation*}
    is homotopy cocartesian.
\end{itemize}
Note that here, $S$ is replaced with $[k] \setminus S$ compared to \autoref{def:new_degree_n_mf_calc} because we consider \emph{covariant} functors instead of contravariant functors.
To conclude, the definition \emph{codegree $n$} of \autoref{def:deg_n_bicartesian} can be seen as a dual to a generalization of the notion of degree $n$ from manifold calculus.

%We can now generalize \autoref{def:new_degree_n_mf_calc} to contravariant functors $F\colon P \to \M$, where $P$ is a lattice and $\M$ is a model category. However, we want to work with covariant functors, not contravariant functors. We therefore replace $S$ with $[k] \setminus S$ and \emph{homotopy cartesian} with \emph{homotopy cocartesian}.
%%\footnote{Making the change from contravariant to covariant functors shouldn't matter. One can translate between covariant and contravariant functors by replacing a poset with its opposite poset, and the opposite of a lattice is still a lattice.}
%\begin{definition}\label{def:degree_n_posets}
%    Let $P$ be a lattice and $\M$ a model category. A functor $F\colon P \to \M$ is \emph{pre-codegree $n$} if, for every $v \in P$ and all $x^0, \dots, x^n$ such that $x^i \vee x^j = v$ for all $i \neq j$, the $n$--cube $\X$ defined by
%    \begin{equation*}
%        \X (S) = 
%        \begin{cases}
%            F(v), &\quad S = [n], \\ 
%            F\left(\bigwedge_{i \notin S} x^i \right), &\quad \text{ otherwise.} \\
%        \end{cases}
%    \end{equation*}
%    is \emph{homotopy cocartesian}.
%\end{definition}
%Note that the fact that $x^i \vee x^j = v$ for $i \neq j$ implies $x^i \le v$ for all $i$.
%
%This definition motivated our definition of degree $n$ functors from lattices (\autoref{def:deg_n_bicartesian}). In particular, it follows from \autoref{prop:defns_consistent} that the two definitions coincide when $P$ is a distributive lattice.

\subsubsection{Notions of dimension}

We will define in \autoref{def:dim} a notion of join-dimension of an element in a lattice. This notion is similar, but not equal, to a notion of dimension introduced in \cite[3.1]{realisationsposetstameness}. The notions differ, for example, in $\Z^2$, where every element has infinite join-dimension in the definition used here, but every element has dimension 2 in the definition used in \cite{realisationsposetstameness}. The two definitions coincide, however, in join-factorization lattices, which we show in \autoref{lemma:dimension_notions_connection}.  

\section{The Taylor telescope for functors on distributive lattices}
\label{sec:taylor_telescope}

In this section, we define the Taylor telescope of a functor and prove \autoref{theoremA_introduction} and \autoref{theoremB_introduction}. We will also give a sufficient condition for the Taylor telescope to converge, and give a brief discussion of poset cocalculus for functors into ordinary categories.

\subsection{Join-decompositions}

An element $x$ in a poset $P$ is \emph{minimal} if there is no other element $y \in P$ such that $y < x$. If $P$ is a lattice, then $x \in P$ is minimal if and only if $x$ is a \emph{least element}, i.e., if $x \le y$ for all $y \in P$. To see this, observe that if $P$ is a lattice and $x \in P$ is minimal, then for every $y \in P$, we have $x \wedge y \le x \Rightarrow x \wedge y = x \Rightarrow x \le y$.
It is obvious that the least element of a poset, if it exists, is unique.

\begin{definition}
    Let $P$ be a lattice, and let $v \in P$. A \emph{join-decomposition} of $v$ (of size $k$) is a finite collection of elements $p^0, \dots, p^{k-1} \in P$ such that $v = p^0 \vee \dots \vee p^{k-1}$.
\end{definition}
%We will sometimes not specify the integer $k$, and just write \emph{decomposition}.
Note that a join-decomposition of size 2 is the same as a pairwise cover of size 2.

Observe that if $p^0, \dots, p^{k-1}$ is a join-decomposition of $v$, we can get a size $k$ pairwise cover of $v$ as follows. Let 
\begin{equation}\label{eq:codecmp_from_decmp}
    x^i = \bigvee_{j \in [k-1] \setminus \{i\}} p^j.
\end{equation}
Then $x^0, \dots, x^{k-1}$ is a pairwise cover of $v$ of size $k$.

\begin{definition}
    Let $x^0, \dots, x^{k-1}$ be a pairwise cover of $v$. We say that the pairwise cover is \emph{reduced} if, for all $i \in [k-1]$, $x^i \neq v$.
\end{definition}

\begin{definition}
    Let $p^0, \dots, p^{k-1}$ be a join-decomposition of $v$. We say that the join-decomposition is \emph{reduced} if none of the $p^i$'s is redundant, i.e., if for all $i$,
    \begin{equation*}
        \bigvee_{j \in [k-1] \setminus \{i\}} p_j \neq v.
    \end{equation*}
    We adopt the convention that $p^0 = v$ is a reduced join-decomposition of $v$, unless $v$ is minimal.
\end{definition}

\begin{example}
Consider the element $(1,1,2,2)$ in the poset $\Nn^4$. A reduced join-decomposition of this element is $\{(1,1,0,0), (1,0,2,1), (0,0,0,2)\}$. The corresponding pairwise cover is $\{(1,1,2,1), (1,1,0,2), (1,0,2,2)\}$.

A join-decomposition of $(1,1,2,2)$ that is not reduced is $\{(1,1,0,0), (0,1,2,0), (0,0,2,2)$, since $(1,1,2,2) = (1,1,0,0) \vee (0,0,2,2)$.
\end{example}

In other words, a join-decomposition is reduced if its corresponding pairwise cover (defined in \eqref{eq:codecmp_from_decmp}) is reduced. If the lattice is distributive, a reduced pairwise cover of size $k$ also gives rise to a reduced join-decomposition. We now show this.

\begin{proposition}
    Let $P$ be a distributive lattice, and let $v \in P$. Let $k$ be a positive integer. Then $v$ has a reduced join-decomposition of size $k$ if and only if it has a reduced pairwise cover of size $k$.
\end{proposition}
\begin{proof}
We first consider the case $k = 1$. An element has a reduced join-decomposition of size 1 if and only if it is not minimal, and it has a reduced pairwise cover of size 1 if and only if it is not minimal, so we are done.

Now, suppose $k \ge 2$. Let $x^0, \dots, x^{k-1}$ be a reduced pairwise cover of $v \in P$. For $i \in [k-1]$, define
\begin{equation}\label{eq:decmp_from_codecmp}
    q^i = \bigwedge_{j \in [k-1] \setminus \{i\}} x^j.
\end{equation}
We show that $q^0, \dots, q^{k-1}$ is a reduced join-decomposition of $v$. To do this, we prove by induction on $n \ge 1$ that
\begin{equation*}
    \bigvee_{i \in \{a_1, \dots, a_n\}} q^i = \X_{x^0, \dots, x^{k-1}}(\{a_1, \dots, a_n\}) =  \bigwedge_{i \notin \{a_1, \dots, a_{n}\}} x^i
\end{equation*}
for any $\{a_1, \dots, a_n\} \subseteq [k-1]$.
This will imply that
\begin{equation*}
    \bigvee_{i \in [k-1]} q^i = v.
\end{equation*}
and
\begin{equation*}
    \bigvee_{i \in [k-1] \setminus \{j\}} q^i = x^j;
\end{equation*}
in other words, that $q^0, \dots, q^{k-1}$ is a reduced join-decomposition of $v$.

The base case in the induction proof follows by definition. The induction step is as follows.
\begin{align*}
    \bigvee_{i \in \{a_1, \dots, a_{n+1}\}} &q^i = q^{a_{n+1}} \vee \bigvee_{i \in \{a_1, \dots, a_{n+1}\}} q^i \\
    &= \left( \bigwedge_{i \notin \{a_{n+1}\}} x^i\right) \vee \left( \bigwedge_{i \notin \{a_1, \dots, a_n\}} x^i \right) \\
    &= \left( \left( \bigwedge_{i \in \{a_1, \dots, a_n\}} x^i\right) \vee x^{a_{n+1}} \right) \wedge \left( \bigwedge_{i \notin \{a_1, \dots, a_{n+1}\}} x^i \right) \quad &\textrm{by \eqref{eq:dist_lattice}}\\
    &= \left( \bigwedge_{i \in \{a_1, \dots, a_n\}} (x^i \vee x^{a_{n+1}}) \right) \wedge \left( \bigwedge_{i \notin \{a_1, \dots, a_{n+1}\}} x^i \right) \quad &\textrm{by \eqref{eq:dist_lattice_alt}}\\
    &= \left( \bigwedge_{i \in \{a_1, \dots, a_n\}} v \right) \wedge \left( \bigwedge_{i \notin \{a_1, \dots, a_{n+1}\}} x^i \right)\\
    &= v \wedge \left( \bigwedge_{i \notin \{a_1, \dots, a_{n+1}\}} x^i \right) = \X_{x^0, \dots, x^{k-1}}(\{a_1, \dots, a_{n+1}\}).\\
\end{align*}
\end{proof}

\begin{remark}
Note that if $q^0, \dots, q^{k-1}$ is a join-decomposition given by \eqref{eq:decmp_from_codecmp} from the pairwise cover $x^0, \dots, x^{k-1}$, and $x^0, \dots, x^{k-1}$ is given by \eqref{eq:codecmp_from_decmp} from the join-decomposition $p^0, \dots, p^{k-1}$, it need not be the case that the two join-decompositions $\{q^0, \dots, q^{k-1}\}$ and $\{p^0, \dots, p^{k-1}\}$ are the same.
In particular, $q^i = \X_{x^0, \dots, x^{k-1}}(\{i\})$ for all $i$, so $q^i \wedge q^j$ has the same value for all $i \neq j$, and this need not be the case for a general join-decomposition.

However, one can show that if it holds that $p^i \wedge p^j$ has the same value for all $i \neq j$, then $q^i = p^i$. % Possible TODO: elaborate, see notes on ipad (remark p,x,q)
\end{remark}

\subsection{Indecomposable join-decompositions}

\begin{definition}
We say that an element $v$ in a lattice $P$ is \emph{join-irreducible} if there don't exist elements $x, y < v$ such that $v = x \vee y$.
\end{definition}

\begin{definition}
We say that a join-decomposition $x^0, \dots, x^k$ is \emph{indecomposable} if it consists only of join-irreducible elements.
\end{definition}

\begin{lemma}{\cite[Lemma 1, p.\ 142]{Birkhoff}}\label{lemma:join_decomp_unique}
    Let $P$ be a distributive lattice. If $v \in P$ has an indecomposable reduced join-decomposition, it is unique.
\end{lemma}
Note that a join-decomposition, by definition, consists of a finite number of elements.

\begin{example}
In the poset $\Nn^4$, the indecomposables are the elements with exactly one nonzero component. The unique indecomposable reduced join-decomposition of $(1, 1, 2, 0) \in \Nn^4$ is $(1,0,0,0), (0,1,0,0), (0,0,2,0)$.
\end{example}

Observe that if an indecomposable join-decomposition is not reduced, it can be reduced by removing redundant elements until it is reduced. Hence, the previous lemma says that if an element has an indecomposable join-decomposition, then it is unique up to adding or removing redundant elements.

The follow definition will be central to the construction of codegree $n$ approximations.
\begin{definition}\label{def:dim}
Let $P$ be a distributive lattice. For an element $v \in P$, we define the \emph{join-dimension} of $v$ as
\begin{equation*}
    \jdim (v) = 
    \begin{cases}
        0, &\ v \textrm{ is minimal}, \\
        k, &\ \text{$v$ is not minimal and has a reduced,}\\[-4pt]
        \ &\ \textrm{indecomposable join-decomposition of size $k$,} \\
        \infty, &\ \text{otherwise.} \\
    \end{cases}
\end{equation*}
\end{definition}
Note that we can consider a least element to be the join of the empty set, and hence it makes sense to set $\jdim(x) = 0$ when $x$ is minimal.
Furthermore, observe that $\jdim(x) \le 1$ if and only if $x$ is join-irreducible.

%\begin{definition}
%    Let $P$ be a lattice. We say that $P$ is a \emph{join-factorization lattice} if it is distributive and every element of $P$ has an indecomposable join-decomposition.
%\end{definition}
%
%Recall the descending chain condition in \autoref{def:descending_chain}. The following proposition shows that there is a large class of distributive lattices that are join-factorization lattices.
%
%\begin{proposition}{\cite[Theorem 9, p.\ 142]{Birkhoff}}\label{thm:Birkhoff}
%    If $P$ be a distributive lattice that satisfies the descending chain condition, then every element of $P$ has a unique indecomposable reduced join-decomposition.
%\end{proposition}

Recall the descending chain condition in \autoref{def:descending_chain}. The following proposition says that in distributive lattices satisfying the descending chain condition, every element has finite join-dimension.

\begin{proposition}{\cite[Theorem 9, p.\ 142]{Birkhoff}}\label{thm:Birkhoff}
    If $P$ is a distributive lattice that satisfies the descending chain condition, then every element of $P$ has a unique indecomposable reduced join-decomposition.
\end{proposition}

\begin{example}\label{ex:factorization_lattices}
Any finite poset satisfies the descending chain condition.
So does a poset $P$ such that every lower level set is finite, i.e., where the set $\{x \in P: x \le v\}$ is finite for any $v \in P$.
Posets of the latter kind include $\Nn, \Nn^k$ and the subposet of finite sets in $\Po(\Nn)$.
\end{example}

\begin{example}
The total order $\R_{\ge 0}$ of nonnegative real numbers does not satisfy the descending chain condition, as the open interval $(0, 1)$ doesn't have a minimal element.
% However, $\R_{\ge 0}$ is a join-factorization lattice.
% Similarly, $(\R_{\ge 0})^k$ for $k \in \Z_{+}$.

% The same applies to $(\R_{\ge 0})^k$ for $k \in \Z_{+}$.
\end{example}

%We will want to work with distributive lattices in which every element has an indecomposable join-decomposition. We give these lattices a name.
We now define a more general class of posets where every element has finite dimension.

\begin{definition}\label{def:join_fact_lattice}
    Let $P$ be a distributive lattice with a minimal element. We say that $P$ is a \emph{join-factorization lattice} if there exists a distributive lattice $Q$, satisfying the descending chain condition, and an order-preserving function $f \colon P \to Q$ such that for all $v \in P$,
    \begin{equation*}
        \jdim(f(v)) = \jdim(v).
    \end{equation*}
\end{definition}

\begin{lemma}
    Let $P$ be a join-factorization lattice. Then every element in $P$ has a unique indecomposable reduced join-decomposition.
\end{lemma}

\begin{proof}
    Uniqueness follows from \autoref{lemma:join_decomp_unique}, as a join-factorization lattice is by definition distributive.
    It now suffices to show that every element in $P$ has finite dimension. Let $x \in P$. By definition, there exists a distributive lattice $Q$ and an order-preserving function $f \colon P \to Q$ such that $\jdim(x) = \jdim(f(x))$. By, \autoref{thm:Birkhoff}, $\jdim(f(x)) < \infty$, which concludes the proof.
\end{proof}

\begin{example}
    Let $P$ be a nonempty distributive lattice satisfying the descending chain condition. Then $P$ is a join-factorization lattice. $P$ has a minimal element, as by the descending chain condition, the subset $P \subseteq P$ must have a minimal element. Thus, $P$ is a join-factorization lattice, as the function $\id_P \colon P \to P$ satisfies the conditions of \autoref{def:join_fact_lattice}.
\end{example}

\begin{example}\label{ex:tot_ord_prod}
% The total order $\R_{\ge 0}$ of nonnegative real numbers does not satisfy the descending chain condition, as the open interval $(0, 1)$ doesn't have a minimal element.
% However, $\R_{\ge 0}$ is a join-factorization lattice.
% Similarly, $(\R_{\ge 0})^k$ for $k \in \Z_{+}$.
% The same applies to $(\R_{\ge 0})^k$ for $k \in \Z_{+}$.
Let $P = T_1 \times \dots \times T_n$, where each $T_i$ is a total order with a minimal element $0$. Then $P$ is a join-factorization lattice. Firstly, $P$ is distributive and $(0, \dots, 0)$ is a minimal element in $P$. Now, consider the order-preserving function $f \colon P \to \{0,1\}^n$ defined by
\begin{equation*}
    f((x_1, \dots, x_n))_i = \begin{cases}
        0, & x_i = 0, \\
        1, & \textrm{otherwise.}
    \end{cases}
\end{equation*}
The join-irreducible elements in $P$ are those with at most one nonzero component, and the join-dimension of an element in $P$ is its number of nonzero components.
Hence, $\jdim(v) = \jdim(f(v))$, and as $\{0, 1 \}^n$ is a finite distributive lattice, this shows that $P$ is a join-factorization lattice.

In particular, $(\nnR)^n$ is a join-factorization lattice.
\end{example}

\begin{example}[Nonexamples]\label{ex:non_factorization_lattices}
The poset $\Po(\Nn)$ doesn't satisfy the descending chain condition. For example, the subset
\begin{equation*}
    \{\{0,1,2,3,\dots\}, \{1,2,3,\dots\}, \{2,3,\dots\}, \dots\}
\end{equation*}
doesn't have a minimal element. Furthermore, $\Po(\Nn)$ is not a join-factorization lattice, as the element $\{0, 1, 2, 3, \dots\}$ can't be written as a finite join of join-irreducible elements (i.e., singletons).

The poset $\Z^2$ is not a join-factorization lattice either, as it doesn't have any join-irreducible elements.

A third distributive lattice that is not a join-factorization lattice is the poset
\begin{equation*}
    P = (\nnR \times \nnR) \setminus \{(0, t) : t > 0\}.
\end{equation*}
Visually, $P$ is the upper right quadrant of the plane, including the x-axis and the origin, but excluding the y-axis (except for the origin). Note that $P$ still has a minimal element, $(0,0)$. The only join-irreducibles in $P$ are of the form $(t, 0)$ for $t \in \R_{\ge 0}$. Hence, the element $(1,1)$, for example, doesn't have an indecomposable join-decomposition.
\end{example}

\begin{lemma}\label{lemma:dimension_notions_connection}
    Let $P$ be a join-factorization lattice and let $v \in P$. Then,    
    \begin{equation*}
        \jdim(v) = \sup \{|U| : U \textrm{ is a reduced join-decomposition of } v\}.
    \end{equation*}
\end{lemma}

\begin{proof}
    Let $x^0, \dots, x^{k-1}$ be a reduced indecomposable join-decomposition of $v$. It suffices to show that any join-decomposition of $v$ with more than $k$ elements is not reduced.

    Let $y^0, \dots, y^{m-1}$ be a join-decomposition of $v$ with $m > k$. Let $i \in \{0, \dots, k-1\}$. Then, as $x^i \le v$,
    \begin{align*}
        x^i &= x^i \wedge v \\
        &= x^i \wedge \bigvee_{j = 0}^{m-1} y^j \\
        &= \bigvee_{j = 0}^{k-1} (x^i \wedge y^j).
    \end{align*}
    Thus, as $x^i$ is join-irreducible, there must exist a $j$ such that $x^i \wedge y^j = x^i$, which implies $x^i \le y^j$.

    Thus, there exists $a_0, \dots, a_{k-1}$ such that $x^i \le y^{a_i}$ for all $i \in \{0, \dots, k-1\}$. Now,
    \begin{align*}
        v = \bigvee_{i=0}^{k-1} x^i
        \le \bigvee_{i=0}^{k-1} y^{a_i} 
        \le \bigvee_{j=0}^{m-1} y^j = v.
    \end{align*}
    Thus,
    \begin{equation*}
        \bigvee_{i=0}^{k-1} y^{a_i} = v,
    \end{equation*}
    which shows that $y^0, \dots, y^{m-1}$ is not redundant.
\end{proof}

\subsection{Codegree $n$ approximations of functors}

We are now ready to define the Taylor telescope of a functor.
Given a distributive lattice $P$, we denote 
\begin{equation}
    P_{\le k} = \{v \in P : \jdim(v) \le k \}.
\end{equation}

\begin{definition}
Let $P$ be a distributive lattice, and $\M$ a good model category. For a functor $F \colon P \to \M$, let $F|_{P_{\le k}}$ denote the functor $F$ restricted to $P_{\le k}$.
We define the \emph{codegree $k$ approximation of $F$}, denoted $T_k F$, as the homotopy left Kan extension of $F|_{P_{\le k}}$ along the inclusion $P_{\le k} \hookrightarrow P$.
On objects, we can write this explicitly as
\begin{equation}\label{eq:dist_lat_TkF}
    T_k F (x) = \underset{v \in P_{\le k}, v \le x}{\hocolim} F(v).
\end{equation}
\end{definition}

From the universal property of $\hocolim$, we get natural transformations $r_k \colon T_k F \to T_{k+1} F$ and $\varepsilon_k \colon T_k F \to F$, such that $\varepsilon_{k+1} \circ r_k = \varepsilon_k$ in $\Ho \Fun(P_{\le k}, \M)$. In total, we get a telescope of functors, that we call the \emph{Taylor telescope} of $F$.
\begin{equation}\label{eq:taylor_tower}
\begin{aligned}
\begin{tikzpicture}
    \diagram{d}{3em}{3em}{
        \vdots & \ \\
        T_1 F & \ \\
        T_0 F & F \\
    };
    
    \path[->,font = \scriptsize, midway]
    (d-2-1) edge node[left]{$r_1$} (d-1-1)
    (d-1-1) edge (d-3-2)
    (d-3-1) edge node[left]{$r_0$} (d-2-1)
    (d-2-1) edge node[below]{$\varepsilon_1$} (d-3-2)
    (d-3-1) edge node[below]{$\varepsilon_0$} (d-3-2);
\end{tikzpicture}
\end{aligned}
\end{equation}
Observe that it follows from \eqref{eq:dist_lat_TkF} that if $x \in P_{\le k}$, then $(\varepsilon_k)_x \colon T_k F (x) \to F(x)$ is a weak equivalence. Furthermore, if $x \in P_{\le k}$, then $(r_k)_x \colon T_k F(x) \to T_{k+1} F(x)$ is a weak equivalence.

We remark that as $\varepsilon_k \colon T_k F \to F$ depends on the functor $F$, and is natural in $F$, a more consistent notation would be to write $(\varepsilon_k)_F \colon T_k F \to F$. However, for notational simplicity, we will usually just write $\varepsilon_k \colon T_k F \to F$.

We will now prove some important properties of the codegree $k$ approximation of a functor. First, we need some lemmas.

\begin{lemma}\label{lemma:irreducible_vee}
    Let $P$ be a distributive lattice, and let $p \in P$ be an join-irreducible element. If $p \le x \vee y$ for some $x, y \in P$, then $p \le x$ or $p \le y$.
\end{lemma}

\begin{proof}
    \begin{equation*}
        p = p \wedge (x \vee y) = (p \wedge x) \vee (p \wedge y).
    \end{equation*}
    Hence, either $p = p \wedge x$ or $p = p \wedge y$, and we are done.
\end{proof}

\begin{lemma}\label{lemma:irreducible_vee_general}
    Let $P$ be a distributive lattice, let $v \in P$ and let $x^0, \dots, x^k$ be a pairwise cover of $v$ of size $k+1$. Let $u \le v$. If $\jdim(u) \le k$, then $u \le x^l$ for some $l \in [k]$.
\end{lemma} 

\begin{proof}
    Let $m = \jdim(u)$, and let $p^0, \dots, p^{m-1}$ be an indecomposable join-decomposition of $u$. Then for all $i \in [m-1]$ and all $0 \le j < l \le k$, we have that $p^i \le x^j \vee x^l$.
    Hence, by \autoref{lemma:irreducible_vee}, for each $i \in [m-1]$, we have that $p^i \le x^j$ for all but one $j \in [k]$.
    Thus, by the pigeonhole principle, there exists an $l \in [k]$ such that $p^i \le x^l$ for all $i \in [m-1]$ (as $m < k+1$). Thus, $u = p^0 \vee \dots \vee p^{k-1} \le x^l$.
\end{proof}

Recall the definition of homotopy finality in \autoref{def:homotopy_finality}. Recall also that a category with an initial object has contractible nerve. In the case where $F \colon C \to D$ is a functor between posets, the under-category $x \downarrow F$ is equivalent to the following sub-poset of $D$.
\begin{equation}\label{eq:under_fiber}
    x \downarrow F \cong \{y \in C : F(y) \ge x\}.
\end{equation}
Recall further that an initial object in a poset is a least element.

\begin{lemma}\label{lemma:final_functor}
    Let $J$ be a lattice, and let $J' \subset J$ be a finite sub-poset. Suppose that $J'$ contains all maximal elements in $J$ and that $J'$ is closed under $\wedge$. Then the inclusion $\iota \colon J' \xhookrightarrow{} J$ is a homotopy final functor.
\end{lemma}

\begin{proof}
We need to show that the poset 
\begin{equation*}
    x \downarrow \iota = \{y \in J' : y \ge x\}
\end{equation*}
has contractible nerve for every $x \in J$. There are two cases to consider.
\begin{description}[style=multiline]
    \item[(I)] $x \in J'$
    \item[(II)] $x \notin J'$
\end{description}
\begin{description}[labelwidth=1.25cm,leftmargin=!]
    \item[Case I] In this case, $x$ is an initial object in $x \downarrow \iota$, and hence $x \downarrow \iota$ has contractible nerve.
    \item[Case II] In this case, observe that $x \downarrow \iota$ is nonempty, as $J'$ contains all maximal elements in $J$. Furthermore, $x \downarrow \iota$ is finite (as $J'$ is finite), so we can define
    \begin{equation*}
        y = \bigwedge_{a \in (x \downarrow \iota)} a
    \end{equation*}
    Then, $y$ is an initial object in $x \downarrow \iota$, and hence $x \downarrow \iota$ has contractible nerve.
\end{description}

\end{proof}

\begin{remark}\label{rmk:rewrite_final_lemma}
The image of a strongly bicartesian cube $\X_{x^0, \dots, x^k} \colon \Po_{k+1} \to P$ is closed under $\wedge$, and so is its punctured version where $\X_{x^0, \dots, x^k}([k])$ is removed. Moreover, $x^0, \dots, x^k$ are the maximal elements in the poset $P' = \{u \in P : u \le x^i \textrm{ for some } i\}$.
Hence, by \autoref{lemma:final_functor} and \autoref{prop:homotopy_finality}, there is a weak equivalence
\begin{equation*}
    \underset{S \subsetneq [k]}{\hocolim} \quad F \circ \X_{x^0, \dots, x^k} (S)
    \simeq \underset{u \leq x^i \textrm{ for some } i}{\hocolim} \quad F (u),
\end{equation*}
for any functor $F \colon P \to \M$.
\end{remark}

\begin{theorem}[Theorem A]\label{theoremA}
    If $P$ is a distributive lattice, then for every functor $F \colon P \to \M$ to a good model category, $T_k F$ is codegree $k$.
\end{theorem}

\begin{proof}
Let $\X$ be a strongly bicartesian $(k+1)$--cube. Using \autoref{lem:codecomp_from_bicartesian}, let $x^0, \dots, x^k$ be the pairwise cover such that $\X = \X_{x^0, \dots, x^k}$. Let $v = \X([k])$.

By \autoref{rmk:rewrite_final_lemma},
\begin{align*}
    &\underset{S \subsetneq [k]}{\hocolim} \quad T_k F \circ \X_{x^0, \dots, x^k} (S) \\
    \simeq &\underset{u \leq x^i \textrm{ for some } i}{\hocolim} \quad T_k F (u). \\
\end{align*}
We can rewrite this as
\begin{align*}
    &\underset{u \leq x^i \textrm{ for some } i}{\hocolim} \quad T_k F (u) \\
    \simeq &\underset{u \leq x^i \textrm{ for some } i}{\hocolim} \quad \underset{w \le u, w \in P_{\le k}}{\hocolim} \quad F(w) \\
    \simeq &\underset{w \leq x^i \textrm{ for some } i, w \in P_{\le k}}{\hocolim} \quad F(w). \\
\end{align*}
By definition, $w \in P_{\le k}$ precisely if $\jdim(w) \le k$. Hence, by \autoref{lemma:irreducible_vee_general}, if $w \in P_{\le k}$ and $w \le v$, then $w \le x^i$ for some $i$. Hence,
\begin{align*}
    &\underset{w \leq x^i \textrm{ for some } i, w \in P_{\le k}}{\hocolim} \quad F(w) \\
    \simeq &\underset{w \leq v, w \in P_{\le k}}{\hocolim} \quad F(w) \\
    \simeq & T_k F (v),
\end{align*}
as desired.

\end{proof}

\begin{theorem}[Theorem B]\label{theoremB}
    Let $P$ be a join-factorization lattice, and let $F \colon P \to \M$ be a functor to a good model category. If $F$ is codegree $n$, then $\varepsilon_n \colon T_n F \to F $ is an objectwise weak equivalence.
\end{theorem}

\begin{proof}
First, let $n = 0$. Let $b$ be the minimal element in $P$. Then, as $F$ is codegree $0$, $F(b \le v)$ is a weak equivalence for all $v \in P$. Hence, as $T_0 F$ is the constant functor at $F(b)$, $\varepsilon_0$ is an objectwise weak equivalence.

Now, let $n \ge 1$.
Let $Q$ be a distributive lattice satisfying the descending chain condition, and $f \colon P \to Q$ an order-preserving function such that $\jdim(v)=\jdim(f(v))$ for all $v \in P$. Let 
\begin{equation*}
    S = \{x \in P : (\varepsilon_n)_x \textrm{ is not a weak equivalence} \}.
\end{equation*}
We prove by contradiction that $S$ is empty.

Suppose that $S$ is nonempty. Let $u$ be a minimal element in $f(S)$ (which exists, as $Q$ satisfies the descending chain condition), and let $v \in f^{-1}(\{u\})$.
Then $v \in S$.
Let further $m = \jdim(v)$.
Then $m > n$, as $(\varepsilon_n)_x$ is a weak equivalence for any $x \in P_{\le n}$.
Now, let $p^0, \dots, p^{m-1}$ be a reduced indecomposable join-decomposition of $v$.
Then $v$ has a pairwise cover $x^0, \dots, x^{m-1}$ given by
\begin{equation*}
    x^i = \bigvee_{j \in [m-1] \setminus \{i\}} p^i.
\end{equation*}
Now, for all $i\in [m-1]$, $\jdim(x^i) = m-1 \neq \jdim(v)$, which implies that $f(x^i) < f(v) = u$. As $u$ is minimal in $S$, $f(x^i) \notin S$, and so $(\varepsilon_n)_{x^i}$ is a weak equivalence. Furthermore, as $f$ is order-preserving, $(\varepsilon_n)_z$ is a weak equivalence for all $z \le x^i$.

Now, as $F$ is codegree $n$, and $m > n$, the cube $F \circ \X_{x^0, \dots, x^{m-1}}$ is homotopy cocartesian.
Furthermore, the cube $T_n F \circ \X_{x^0, \dots, x^{m-1}}$ is homotopy cocartesian.
Thus, as $(\varepsilon_n)_z$ is an objectwise weak equivalence for all
\begin{equation*}
    z \in \{\X_{x^0, \dots, x^{m-1}}(S) : S \subsetneq [m-1]\},
\end{equation*}
the morphism $(\varepsilon_n)_v$ is a weak equivalence. % (as $v = \X_{x^0, \dots, x^{m-1}}([m-1]))$.
This contradicts the fact that $v \in S$.
Hence, $S$ is empty, which completes the proof.

\end{proof}

The following proposition states that when $P$ is a join-factorization lattice, then the Taylor telescope of $F$ converges to $F$.

\begin{proposition}
    Let $P$ be a join-factorization lattice, and let $F \colon P \to \M$ be a functor to a good model category.
    For $x \in P$, let $\hocolim_k T_k F (x)$ denote the homotopy colimit of the telescope $T_0 F (x) \to T_1 F (x) \to T_2 F (x) \to \cdots$.
    The induced map
    \begin{equation*}
        \hocolim_k T_k F (x) \to F (x)
    \end{equation*}
    is a weak equivalence.
\end{proposition}

\begin{proof}
The fact that $P$ is a join-factorization lattice implies that every element of $P$ has finite join-dimension. When $\jdim(x) = N$, then $(\varepsilon_m)_x \colon T_m F (x) \to F (x)$ are weak equivalences for every $m \ge N$. The result now follows from \autoref{lemma:telescope_hocolim}.
\end{proof}

Finally, the following proposition shows that the functor $T_n F$ is in some sense the \emph{universal} codegree $n$ approximation of a functor $F$.

\begin{proposition}
    Let $P$ be a join-factorization lattice, $\M$ a good model category, and $F \colon P \to \M$ a functor.
    Then $\varepsilon_n \colon T_n F \to F$ is the universal map from a codegree $n$ functor to $F$, in the homotopy category $\Ho \Fun(P, \M)$.

    In other words, for any codegree $n$ functor $G \colon P \to \M$ and natural transformation $\zeta \colon G \to F$, there exists a unique natural transformation $\upsilon \colon G \to T_n F$ in $\Ho \Fun(P, \M)$ such that $\zeta = \varepsilon_n \circ \upsilon$.
\end{proposition}

The proof resembles that of \cite[Theorem 1.8]{goodwillie_calculus3}.

\begin{proof}
    Let $F \colon P \to \M$, and let
    $\zeta \colon G \to F$ be a natural transformation from a codegree $n$ functor to $F$.
    Consider the commutative diagram
    \begin{center}
    \begin{tikzcd}
    T_n G \arrow[r, "T_n \zeta"] \arrow[d, "(\varepsilon_n)_G"] & T_n F \arrow[d, "(\varepsilon_n)_F"] \\
    G \arrow[r, "\zeta"]                                 & F                            
    \end{tikzcd}
    \end{center}
    As $(\varepsilon_n)_G$ is a weak equivalence, it is an isomorphism in the homotopy category $\Ho \Fun(P, \M)$. Hence, we can take $\upsilon = T_n \zeta \circ ((\varepsilon_n)_G)^{-1}$, which satisfies $(\varepsilon_n)_F \circ \upsilon = \zeta$.
    
    It remains to show uniqueness. Suppose that $(\varepsilon_n)_F \circ \upsilon = \zeta$ and $(\varepsilon_n)_F \circ \upsilon' = \zeta$. Then $T_n \zeta = T_n ((\varepsilon_n)_F \circ \upsilon) = T_n (\varepsilon_n)_F \circ T_n \upsilon$, and likewise, $T_n \zeta = T_n (\varepsilon_n)_F \circ T_n \upsilon'$. As $T_n (\varepsilon_n)_F$ is an isomorphism in the homotopy category, $T_n \upsilon = T_n \upsilon'$. Thus, as $\upsilon$ and $\upsilon'$ are natural transformations between codegree n functors, $\upsilon = \upsilon'$.
\end{proof}

\subsection{Poset cocalculus for ordinary categories}

Note that for any bicomplete category $\C$, we can apply poset cocalculus to functors from a poset into $\C$, by considering the trivial model structure on $\C$ where the weak equivalences are the isomorphisms in $\C$. In this case, $\hocolim$ of any diagram is just $\colim$.
Furthermore, it is easy to verify that it suffices to require that $\C$ is \emph{cocomplete}, instead of bicomplete, for the relevant definitions and theorems to hold.
We recall the most important parts of poset cocalculus, rephrased in this specific setting.

\begin{definition}\label{def:deg_n_bicartesian_cat}
    Let $P$ be a lattice and $\C$ a cocomplete category. A functor $F \colon P \to \C$ is \emph{codegree $n$} if it sends strongly bicartesian $(n+1)$--cubes to cocartesian cubes.
\end{definition}

The codegree $k$ approximation of a functor $F \colon P \to \C$, is defined as
\begin{equation}\label{eq:dist_lat_TkF_cat}
    T_k F (x) = \underset{v \in P_{\le k}, v \le x}{\colim} F(v).
\end{equation}

\begin{theorem}[Theorem A for categories]\label{theoremA_cat}
    If $P$ is a distributive lattice, then for every functor $F \colon P \to \C$ to a cocomplete category, $T_k F$ is codegree $k$.
\end{theorem}

\begin{theorem}[Theorem B for categories]\label{theoremB_cat}
    Let $P$ be a join-factorization lattice, and let $F \colon P \to \C$ be a functor to a cocomplete category. If $F$ is codegree $n$, then $\varepsilon_n \colon T_n F \to F $ is a natural isomorphism.
\end{theorem}

\begin{remark}
Let $P$ be a join-factorization lattice and $\C$ a cocomplete category. Let $\Fun_n(P, \C)$ denote the full subcategory of $\Fun(P, \C)$ consisting of codegree $n$ functors.
Then $T_n$ is the right adjoint of the inclusion $\iota \colon \Fun_n(P, \C) \to \Fun(P,\C)$, i.e., we have an adjunction
\[
    \begin{tikzcd}[column sep=large]
        \Fun_n(P, \C) \arrow[r, shift left=1ex, "\iota"{name=G, yshift=1pt}] & \Fun(P, \C) \arrow[l, shift left=.5ex, "T_n"{name=F}]
        \arrow[phantom, from=F, to=G, , "\scriptscriptstyle\boldsymbol{\bot}"].
    \end{tikzcd}
\]
The counit 
\begin{equation*}
    \hat \varepsilon \colon \iota \circ T_n \to \id_{\Fun(P,\C)}
\end{equation*}
is defined objectwise as the natural transformation
\begin{equation*}
    \hat \varepsilon_F = \varepsilon_n \colon T_n F \to F.
\end{equation*}
Recall that $\varepsilon_n$ is an isomorphism whenever $F \in \Fun_n(P, \C)$, by \autoref{theoremB_cat}.
The unit 
\begin{equation*}
    \hat \eta \colon \id_{\Fun_n(P, \C)} \to T_n \circ \iota
\end{equation*}
is defined objectwise as
\begin{equation*}
    \hat \eta_F = (\varepsilon_n)^{-1} \colon F \to T_n F.
\end{equation*}
The naturality of the unit and counit follows from the universal property of left Kan extensions. The counit-unit equations follow from the fact that $\hat \varepsilon_{T_n F} = T_n \hat \varepsilon_F$, which also follow from the universal property of left Kan extensions.
% TODO: should I spell this out more? It is quite interesting
\end{remark}

\section{Examples of the codegree $n$ approximation}

We give examples of the codegree $n$ approximation of a functor.

\subsection{Examples}
We consider examples with lattices that are finite products of total orders with minimal elements, i.e., lattices of the form $P = P_1 \times \cdots \times P_n$, where each $P_i$ is a total order with a minimal element 0. These lattices are join-factorization lattices by \autoref{ex:tot_ord_prod}, and the indecomposable join-decomposition of an element $(v_1, \dots, v_n)$ is
\begin{equation*}
    \{(v_1, 0, \dots, 0), (0, v_2, 0, \dots, 0), \dots, (0, \dots, 0, v_n)\}. 
\end{equation*}
The join-dimension of an element $(v_1, \dots, v_n)$ is the number of nonzero components.

\begin{remark}\label{remark_lemma1}
\autoref{lemma:final_functor} implies that for $P = P_1 \times \dots \times P_n$ a product of total orders with minimal elements and $F \colon P \to \M$ a functor into a good model category, we can compute the homotopy colimit $T_k F (x)$ from just the elements of the form $(v_1, \dots, v_n)$ where each $v_i$ is equal to either $x_i$ or 0. We can therefore give an equivalent formula as follows.

For a set $S \subset \{1, \dots, m\}$ and $x \in P$, let $\lambda(x, S)$ be the element in $P$ given by
\begin{equation*}
    (\lambda(x, S))_i = 
    \begin{cases}
        x_i, &\quad i \in S, \\ 
        0, &\quad i \notin S.\\ 
    \end{cases}
\end{equation*}
For example, $\lambda ((x_1, x_2, x_3), \{1, 3\}) = (x_1, 0, x_3)$.
Then,
\begin{equation}\label{eq:alt_def_taylor_tower_posets}
    T_k F(x) \simeq \underset{S \subset \{1, \dots, m\}, |S| \le k}{\hocolim} F(\lambda (x, S)).
\end{equation}
In particular, $T_0 F$ is the constant functor at  $F(0, \dots, 0)$, i.e., $T_0 F(x) \simeq F(0, \dots, 0)$ for all $x \in P$.
\end{remark}

% Example: some \R^2 stuff
\begin{example}
For the poset $P = \nnR \times \nnR$, we have
\begin{equation*}
    P_{\le k} =
    \begin{cases}
        0, &\quad k = 0, \\ 
        (\nnR \times \{0\}) \cup (\{0\} \times \nnR), &\quad k = 1, \\ 
        P, &\quad k \ge 2, \\ 
    \end{cases}
\end{equation*}
In general, for $P = (\nnR)^{\times n}$, we have that $P_{\le k}$ is the union of $n \choose k$ orthogonal $k$--planes in $\R^n$ (restricted to $(\nnR)^{\times n}$).
% Illustrations would be nice here
\end{example}

\begin{example}\label{ex:simple_t1_example}
We compute $T_1$ of the functor in \autoref{ex:deg_2_functor}. The join-dimension 1 elements are $(0, 0), (0,1)$ and $(1,0)$. Thus, $T_1 F (v) \simeq F(v)$ everywhere except in $v = (1,1)$. We compute $T_1 F((1,1))$ as the homotopy colimit of:
\begin{center}
\begin{tikzpicture}
    \diagram{d}{3em}{3em}{
        \Fb & \  \\
        0 & \Fb \\
    };
    
    \path[->,font = \scriptsize, midway]
    (d-2-1) edge (d-1-1)
    (d-2-1) edge (d-2-2);
\end{tikzpicture}
\end{center}
We get that $T_1 F((1,1)) = \Fb \oplus \Fb$, and $T_1 F$ is:
\begin{center}
\begin{tikzpicture}
    \diagram{d}{3em}{3em}{
        \Fb & \Fb \oplus \Fb \\
        0 & \Fb \\
    };
    
    \path[->,font = \scriptsize, midway]
    (d-1-1) edge node[above]{$\begin{pmatrix}1 \\ 0\end{pmatrix}$} (d-1-2)
    (d-2-2) edge node[right]{$\begin{pmatrix}0 \\ 1\end{pmatrix}$} (d-1-2)
    (d-2-1) edge (d-1-1)
    (d-2-1) edge (d-2-2);
\end{tikzpicture}
\end{center}

\end{example}

\begin{example}
We look at functors $\{0, 1\}^2 \to \Ch_{\Fb}$. First, let $S^n$ denote the chain complex with a single copy of $\Fb$ in degree $n$ and zero everywhere else. Let further $D^n$ denote the chain complex with $\Fb$ in degree $n$ and degree $n-1$ and zero everywhere else, and where the differential in degree $n$ is the identity.

We compute $T_1$ of the following functor, $F \colon \{0,1\}^2 \to \Ch_{\Fb}$,
\begin{center}
\begin{tikzpicture}
    \diagram{d}{3em}{3em}{
        0 & 0 \\
        S^0 & 0. \\
    };
    
    \path[->,font = \scriptsize, midway]
    (d-1-1) edge (d-1-2)
    (d-2-2) edge (d-1-2)
    (d-2-1) edge (d-1-1)
    (d-2-1) edge (d-2-2);
\end{tikzpicture}
\end{center}

As in \autoref{ex:simple_t1_example}, we need to compute $T_1 F(1,1)$, which is the homotopy colimit of
\begin{center}
\begin{tikzpicture}
    \diagram{d}{3em}{3em}{
        0 & \ \\
        S^0 & 0. \\
    };
    
    \path[->,font = \scriptsize, midway]
    (d-2-1) edge (d-1-1)
    (d-2-1) edge (d-2-2);
\end{tikzpicture}
\end{center}
We compute a cofibrant replacement of this diagram,
\begin{center}
\begin{tikzpicture}
    \diagram{d}{3em}{3em}{
        D^1 & \ \\
        S^0 & D^1. \\
    };
    
    \path[>->,font = \scriptsize, midway]
    (d-2-1) edge (d-1-1)
    (d-2-1) edge (d-2-2);
\end{tikzpicture}
\end{center}
The colimit of this diagram is the chain complex
\begin{center}
\begin{tikzpicture}
    \diagram{d}{2em}{2em}{
        \Fb \oplus \Fb \\
        \Fb, \\
    };
    
    \path[->,font = \scriptsize, midway]
    (d-1-1) edge node[left]{$(1 \  1)$} (d-2-1);
\end{tikzpicture}
\end{center}
which is weakly equivalent to $S^1$ through the chain map
\begin{center}
\begin{tikzpicture}
    \diagram{d}{2em}{4em}{
        \Fb \oplus \Fb & \Fb \\
        \Fb & 0. \\
    };
    
    \path[->,font = \scriptsize, midway]
    (d-1-1) edge node[left]{$(1 \ 1)$} (d-2-1)
    (d-1-2) edge node[right]{$0$} (d-2-2);
    
    \path[->,font = \scriptsize, midway, dashed]
    (d-1-1) edge node[above]{$(1 \ 1)$} (d-1-2)
    (d-2-1) edge node[above]{$0$} (d-2-2);
\end{tikzpicture}
\end{center}

Hence, $T_1 F$ is
\begin{center}
\begin{tikzpicture}
    \diagram{d}{3em}{3em}{
        D^1 & S^1 \\
        S^0 & D^1. \\
    };
    
    \path[->,font = \scriptsize, midway]
    (d-1-1) edge (d-1-2)
    (d-2-2) edge (d-1-2)
    (d-2-1) edge (d-1-1)
    (d-2-1) edge (d-2-2);
\end{tikzpicture}
\end{center}

We can verify that this functor is codegree 1. Indeed, $S^1$ is weakly equivalent to the mapping cone of $S^0 \to D^1 \oplus D^1$, which is what we need to check.
\end{example}

\begin{example}\label{ex:nonpolynomial_functor}
We give an example of a non-polynomial functor, i.e., a functor that is not codegree $n$ for any $n \in \Z_+$.

Let $P$ be the poset of finite subsets of $\Nn$, ordered by inclusion. Then $P$ is a join-factorization lattice, as noted in \autoref{ex:factorization_lattices}.
Let further $[0, \infty]$ be the poset of nonnegative real numbers adjoined with a greatest element $\infty$. The poset $[0, \infty]$ admits all joins and meets, also infinite ones, and is thus a bicomplete category. Let 
\begin{equation*}
    F \colon P \to [0, \infty]
\end{equation*}
be defined by
\begin{equation*}
    F(X) = |X|.
\end{equation*}

Then,
\begin{align*}
    T_m F (\{1, \dots, n\}) &= \sup_{X \subseteq \{1, \dots, n\}, |X| \le m} |X| \\
    &=
    \begin{cases}
        n, &\quad n \le m,  \\ 
        m, &\quad \textrm{otherwise.} \\ 
    \end{cases}
\end{align*}
In particular, for all $n \in \Z_+$, we have that $T_n F(\{1, \dots, n+1\}) \neq F(\{1, \dots, n + 1\})$. Thus, by \autoref{theoremB_cat}, $F$ is not codegree $n$ for any $n$.
\end{example}

\subsection{Nonexamples}

\begin{example}
\autoref{fig:M3_and_N5} shows two non-distributive lattices, $M3$ and $N5$. A lattice is distributive if and only if it doesn't contain $M3$ or $N5$ as a sublattice \cite[Chapter IX, Theorem 2]{Birkhoff}. We give an example where \autoref{theoremA} fails to hold for each of these lattices.\footnote{Note that $T_k$ of a functor $F \colon P \to \M$ is only defined when $P$ is a distributive lattice. However, $T_1$ is easily generalized to arbitrary lattices, by letting $P_{\le 1}$ be the set of join-irreducible elements.}

\begin{figure}[htbp]
    \centering
    \begin{subfigure}{0.45\textwidth}
        \centering
        \begin{tikzpicture}
            \diagram{d}{1.5em}{3em}{
                \ & a_4 & \ \\
                \ & \ & \ \\
                a_1 & a_2 & a_3 \\
                \ & \ & \ \\
                \ & a_0 & \ \\
            };
            \path[-,font = \scriptsize, midway]
            (d-3-1) edge (d-1-2)
            (d-3-2) edge (d-1-2)
            (d-3-3) edge (d-1-2)
            (d-3-1) edge (d-5-2)
            (d-3-2) edge (d-5-2)
            (d-3-3) edge (d-5-2);
        \end{tikzpicture}
        \caption{$M_3$}
    \end{subfigure}
    \hfill
    \begin{subfigure}{0.45\textwidth}
        \centering
        \begin{tikzpicture}
            \diagram{d}{1.5em}{3em}{
                \ & b_4 & \ \\
                b_3 & \ & \ \\
                \ & \ & b_2 \\
                b_1 & \ & \ \\
                \ & b_0 & \ \\
            };
            
            \path[-,font = \scriptsize, midway]
            (d-2-1) edge (d-1-2)
            (d-4-1) edge (d-2-1)
            (d-3-3) edge (d-1-2)
            (d-5-2) edge (d-4-1)
            (d-5-2) edge (d-3-3);
        \end{tikzpicture}
        \caption{$N_5$}
    \end{subfigure}
    \caption{The Hasse diagrams of two non-distributive lattices.}
    \label{fig:M3_and_N5}
\end{figure}
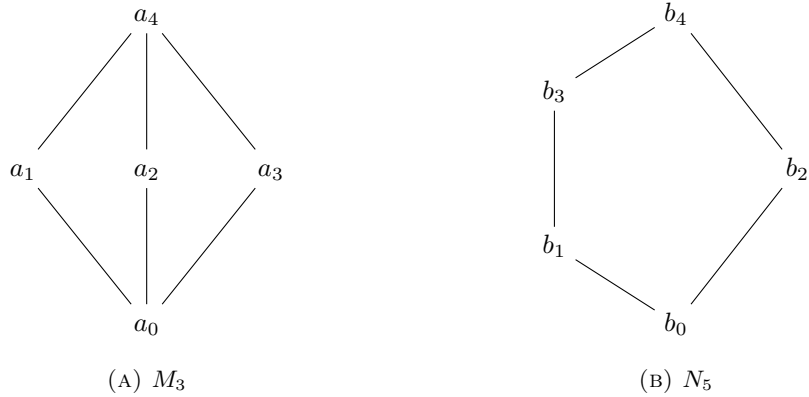

We first consider the following functor $F \colon M_3 \to \Vecs_{\Fb}$.
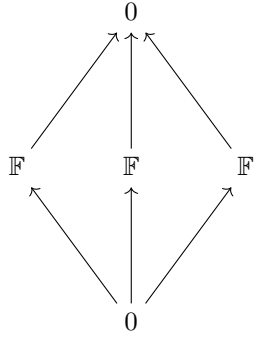
\begin{center}
\begin{tikzpicture}
    \diagram{d}{1.5em}{3em}{
        \ & 0 & \ \\
        \ & \ & \ \\
        \Fb & \Fb & \Fb \\
        \ & \ & \ \\
        \ & 0 & \ \\
    };
    \path[->,font = \scriptsize, midway]
    (d-3-1) edge (d-1-2)
    (d-3-2) edge (d-1-2)
    (d-3-3) edge (d-1-2)
    (d-5-2) edge (d-3-1)
    (d-5-2) edge (d-3-2)
    (d-5-2) edge (d-3-3);
\end{tikzpicture}
\end{center}
We have that $(M_3)_{\le 1} = \{a_0,a_1,a_2,a_3\}$, so $T_1 F(x) \cong F(x)$ when $x \neq a_4$. Furthermore, 
\begin{align*}
    T_1 F(a_4) &= \underset{x \le a_4, x \in (M_3)_{\le 1}}{\colim} F(x) \\
    &= \underset{x \in \{a_0, a_1, a_2, a_3\}}{\colim} F(x) = \Fb^3.
\end{align*}
Thus, $T_1 F$ is
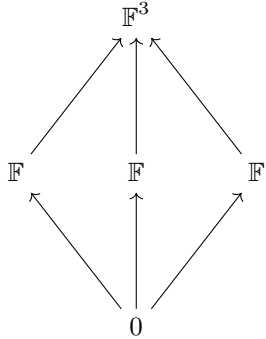
\begin{center}
\begin{tikzpicture}
    \diagram{d}{1.5em}{3em}{
        \ & \Fb^3 & \ \\
        \ & \ & \ \\
        \Fb & \Fb & \Fb \\
        \ & \ & \ \\
        \ & 0 & \ \\
    };
    \path[->,font = \scriptsize, midway]
    (d-3-1) edge (d-1-2)
    (d-3-2) edge (d-1-2)
    (d-3-3) edge (d-1-2)
    (d-5-2) edge (d-3-1)
    (d-5-2) edge (d-3-2)
    (d-5-2) edge (d-3-3);
\end{tikzpicture}
\end{center}
We show that the 2--cube $(T_1 F) \circ \X_{a_1, a_2}$ is not cocartesian, which will show that $T_1 F$ is not cocartesian.
The cube $\X_{a_1, a_2}$ consists of $a_1$, $a_2$, $a_1 \wedge a_2 = a_0$ and $a_1 \vee a_2 = a_4$. Thus, $(T_1 F) \circ \X_{a_1, a_2}$ is
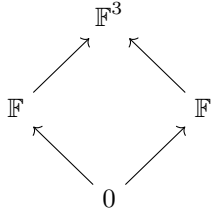
\begin{center}
\begin{tikzpicture}
    \diagram{d}{2em}{2em}{
        \ & \Fb^3 & \ \\
        \Fb & \ & \Fb \\
        \ & 0 & \ \\
    };
    \path[->,font = \scriptsize, midway]
    (d-2-1) edge (d-1-2)
    (d-2-3) edge (d-1-2)
    (d-3-2) edge (d-2-1)
    (d-3-2) edge (d-2-3);
\end{tikzpicture}
\end{center}
As $\Fb^3$ is not isomorphic to $\Fb^2$, this is not cocartesian.

Next, we consider the following functor $G \colon N_5 \to \Vecs_{\Fb}$
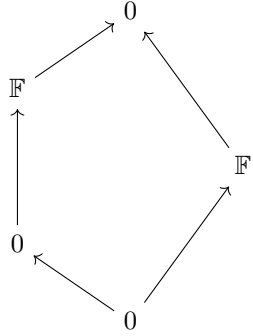
\begin{center}
\begin{tikzpicture}
    \diagram{d}{1.5em}{3em}{
        \ & 0 & \ \\
        \Fb & \ & \ \\
        \ & \ & \Fb \\
        0 & \ & \ \\
        \ & 0 & \ \\
    };
    
    \path[->,font = \scriptsize, midway]
    (d-2-1) edge (d-1-2)
    (d-4-1) edge (d-2-1)
    (d-3-3) edge (d-1-2)
    (d-5-2) edge (d-4-1)
    (d-5-2) edge (d-3-3);
\end{tikzpicture}
\end{center}
We have that $(N_5)_{\le 1} = \{b_0, b_1, b_2, b_3\}$. By similar computations as before, we get that $T_1 G$ is
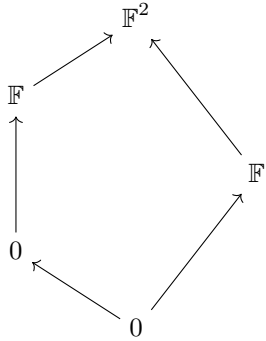
\begin{center}
\begin{tikzpicture}
    \diagram{d}{1.5em}{3em}{
        \ & \Fb^2 & \ \\
        \Fb & \ & \ \\
        \ & \ & \Fb \\
        0 & \ & \ \\
        \ & 0 & \ \\
    };
    
    \path[->,font = \scriptsize, midway]
    (d-2-1) edge (d-1-2)
    (d-4-1) edge (d-2-1)
    (d-3-3) edge (d-1-2)
    (d-5-2) edge (d-4-1)
    (d-5-2) edge (d-3-3);
\end{tikzpicture}
\end{center}
This is again not codegree 1, as the 2--cube $(T_1 G) \circ \X_{b_1, b_2}$ is not cocartesian.
\end{example}

The following two examples exhibit functors on distributive lattices that are not join-factorization lattices, and where \autoref{theoremB} fails to hold.

\begin{example}
Let $P$ be the poset
\begin{equation*}
    P = (\nnR \times \nnR) \setminus \{(0, t) : t > 0\}
\end{equation*}
from \autoref{ex:non_factorization_lattices}.
This is a distributive lattices, but not a join-factorization lattice, and the join-irreducible elements are those on the form $(t, 0), t \in \R_{\ge 0}$.
Hence, the elements that have a positive value in both coordinates don't have an indecomposable join-decomposition.
Thus, we have $P_{\le 0} = \{(0,0)\}$ and $P_{\le 1} = \{(t, 0) : t \ge 0\}$, and all other elements have infinite join-dimension.

Now consider the functor $F \colon P \to \Vecs_{\Fb}$ defined by 
\begin{equation*}
    F((x,y)) =
    \begin{cases}
        \Fb, &\quad y=1, \\ 
        0, &\quad \textrm{otherwise.} \\ 
    \end{cases}
\end{equation*}
We have that $F$ is codegree 1, but $T_1 F = 0$ everywhere.
\end{example}

\begin{example}\label{ex:infinite_product}
We show that the telescope in \eqref{eq:taylor_tower} need not converge for functors on distributive lattices that are not join-factorization lattices. Consider the poset $\Po(\Nn)$ of (not necessarily finite) subsets of $\Nn$, ordered by inclusion. In this poset, infinite subsets have infinite join-dimension. Consider a functor $F \colon \Po(\Nn) \to \Vecs_{\Fb}$ given by
\begin{equation*}
    F(X) =
    \begin{cases}
        0, &\quad \textrm{$X$ is finite,} \\ 
        \Fb, &\quad \textrm{$X$ is infinite.} \\ 
    \end{cases}
\end{equation*}
As $F(x) = 0$ for all $x$ with finite join-dimension, we see that $T_m F = 0$ for all $m$. Hence, the telescope clearly doesn't converge to $F$.
\end{example}

\section{Stability under the interleaving distance}

%\subsection{Stability}

In this section, we will consider functors $F \colon (\nnR)^n \to \M$, where $\M$ is a good model category. The results here thus applies both to functors into $\Ch_{\Fb}$ (equipped with the injective model structure) and into $\Vecs_{\Fb}$ (equipped with the trivial model structure).

\subsection{Generalized homotopy interleaving distance}

We combine here the notion of generalized interleaving distance from \cite{generalizedInterleaving} with the concept of homotopy interleaving distance from \cite{blumberg2022universality}.

A \emph{translation} on a poset $P$ is a functor $\Gamma \colon P \to P$ satisfying $x \le \Gamma(x)$. Let $\Trans_P$ be the set of translations on $P$. This is a poset under the relation
\begin{equation*}
    \Gamma \le K \Leftrightarrow \Gamma(x) \le K(x) \textrm{ for all $x$.} 
\end{equation*}
We let $\eta_{\Gamma} \colon \id_P \to \Gamma$ denote the natural transformation given by $(\eta_{\Gamma})_x \colon x \le \Gamma(x)$.

Let $\Gamma \in \Trans_P$, and let $F, G$ be functors from $P$ to a category $\C$. An \emph{$\Gamma$--interleaving} is a pair of natural transformations
\begin{equation*}
    \phi \colon F \to G \circ \Gamma, \quad
    \psi \colon G \to F \circ \Gamma,
\end{equation*}
such that $\psi_{\Gamma} \circ \phi = F (\eta_{\Gamma^2})$ and $\phi_{\Gamma} \circ \psi = G (\eta_{\Gamma^2})$, where $\Gamma^2 = \Gamma \circ \Gamma$.

A \emph{superlinear family} is a function $\Omega \colon [0, \infty) \to \Trans_P$ that satisfies $\Omega_{\varepsilon_1 + \varepsilon_2} \ge \Omega_{\varepsilon_2} \Omega_{\varepsilon_2}$ for all $\varepsilon_1, \varepsilon_2 \ge 0$. Given a superlinear family $\Omega \colon [0, \infty) \to \Trans_P$, we define the \emph{interleaving distance} between two functors $F, G \colon P \to \C$ as
\begin{equation*}
    d^{\Omega}(F, G) = \inf\{\varepsilon \in [0, \infty) : F, G \textrm{ are $\Omega_\varepsilon$--interleaved}\}.
\end{equation*}
This is a pseudometric by \cite[Theorem 2.5.3]{generalizedInterleaving}.

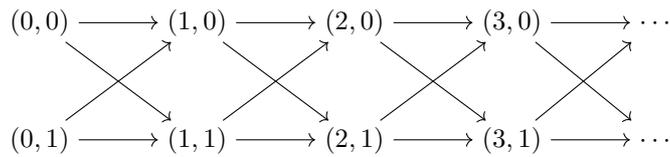
\begin{figure}[hbtp]
    \centering
    \begin{tikzpicture}
        \diagram{d}{3em}{3em}{
            (0,0) & (1,0) & (2,0) & (3,0) & \cdots \\
            (0,1) & (1,1) & (2,1) & (3,1) & \cdots \\
        };
        
        \path[->,font = \scriptsize, midway]
        (d-1-1) edge (d-1-2)
        (d-1-2) edge (d-1-3)
        (d-1-3) edge (d-1-4)
        (d-1-4) edge (d-1-5)
        (d-2-1) edge (d-2-2)
        (d-2-2) edge (d-2-3)
        (d-2-3) edge (d-2-4)
        (d-2-4) edge (d-2-5)
        (d-1-1) edge (d-2-2)
        (d-1-2) edge (d-2-3)
        (d-1-3) edge (d-2-4)
        (d-1-4) edge (d-2-5)
        (d-2-1) edge (d-1-2)
        (d-2-2) edge (d-1-3)
        (d-2-3) edge (d-1-4)
        (d-2-4) edge (d-1-5);
    \end{tikzpicture}
    \caption{An illustration of the poset $\Ib^{\Gamma}$, where $P = \Nn$, and $\Gamma \in \Trans_P$ is the translation given by $\Gamma(n) = n+1$.}
    \label{fig:interleaving_poset}
\end{figure}

A $\Gamma$--interleaving can equivalently be formulated as follows. Given $\Gamma \in \Trans_P$, let $\Ib^{\Gamma}$ denote the poset with object set $P \times \{0, 1\}$ and poset relation: $(x, i) \le (y, j)$ if either
\begin{itemize}
    \item $i = j$ and $x \le y$, or
    \item $\Gamma(x) \le y$.
\end{itemize}
Let $E^0 \colon P \to \Ib^{\Gamma}$ and $E^1 \colon P \to \Ib^{\Gamma}$ be the inclusions into $P \times \{0\}$ and $P \times \{1\}$, respectively. Then a $\Gamma$--interleaving between $F, G \colon P \to \C$ is a functor $Z \colon \Ib^\Gamma \to \C$ such that $Z \circ E^0 = F$ and $Z \circ E^1 = G$.
An illustration of the poset $\Ib^{\Gamma}$ is given in \autoref{fig:interleaving_poset}.

Finally, we define the notion of \emph{homotopy interleaving}. For two objects $X, Y$ in a category with weak equivalences, we write $X \simeq Y$ if $X$ and $Y$ are connected by a zigzag of weak equivalences.
\begin{definition}
    Let $F, G \colon P \to \M$ be functors from a poset to a model category, and let $\Gamma \in \Trans_P$. We say that $F$ and $G$ are \emph{$\Gamma$--homotopy-interleaved} if there exists $F' \simeq F$ and $G' \simeq G$ such that $F'$ and $G'$ are $\Gamma$--interleaved.
\end{definition}
Given a superlinear family $\Omega \colon [0, \infty) \to \Trans_P$, we define the \emph{$\Omega$--homotopy interleaving distance} between two functors $F, G \colon P \to \C$, where $\C$ is a model category, as
\begin{equation*}
    d^{\Omega}_H (F, G) = \inf\{\varepsilon \in [0, \infty) : F, G \textrm{ are $\Omega_\varepsilon$--homotopy-interleaved}\}.
\end{equation*}

\subsection{Multiplicative interleaving distance}

We now restrict ourselves to functors from $(\R_{\ge 0})^n$. % maybe TODO: give example of why our approximations are not stable with respect to ordinary interleaving distance

\begin{definition}
    We define $\Lambda^{(n)} \colon [0, \infty) \to \Trans_{(\R_{\ge 0})^n}$ as the superlinear family given by
    \begin{align*}
        \Lambda^{(n)}_{\varepsilon} \colon \qquad (\R_{\ge 0})^n &\to (\R_{\ge 0})^n \\
        (v_1, \dots, v_n) &\mapsto (v_1 \cdot e^{\varepsilon}, \dots, v_n \cdot e^{\varepsilon}).
    \end{align*}
    When it's clear from the context what $n$ is, we will omit the upper index and just write $\Lambda$.
\end{definition}

A direct computation shows that $\Lambda_{\varepsilon_1} \Lambda_{\varepsilon_2} = \Lambda_{\varepsilon_1 + \varepsilon_2}$, so that $\Lambda$ is a superlinear family.

\begin{definition}
    Let $F, G \colon (\R_{\ge 0})^n \to \M$, where $\M$ is a model category. We define the \emph{multiplicative interleaving distance} between $F$ and $G$ as
    \begin{equation*}
        d^{\Lambda}_H(F, G) = \inf \{\varepsilon \in [0, \infty) : \textrm{$F$ and $G$ are $\Lambda_{\varepsilon}$--homotopy-interleaved} \}.
    \end{equation*}
\end{definition}

\subsection{Stability of codegree $k$ approximations}

Observe that $\Lv(P_{\le k}) \subseteq P_{\le k}$ (where $P = (\R_{\ge 0})^n$). This key property allows us to prove the following proposition.

Recall that given $F \colon I \to \C$ and $\alpha \colon I \to J$, $\Lan_{\alpha} F$ denotes the left Kan extension of $F$ along $\alpha$.

\begin{proposition}\label{prop:interleaving_lke}
    Let $P = (\R_{\ge 0})^n$, and let $F, G \colon P \to \C$ be functors to a cocomplete category. Suppose that the pair $(\phi, \psi)$, where
    \begin{align*}
        \phi \colon F &\to G \circ \Lv,\\
        \psi \colon G &\to F \circ \Lv,
    \end{align*}
    is a $\Lv$--interleaving between $F$ and $G$. Then there is an induced $\Lv$--interleaving $(\hat \phi, \hat \psi)$, where
    \begin{align*}
        \hat \phi \colon \Lan_i (F|_{P_{\le k}}) &\to (\Lan_i (G|_{P_{\le k}})) \circ \Lv,\\
        \hat \psi \colon \Lan_i (G|_{P_{\le k}}) &\to (\Lan_i (F|_{P_{\le k}})) \circ \Lv,
    \end{align*}
    between $\Lan_i(F|_{P_{\le k}})$ and $\Lan_i (G|_{P_{\le k}})$, where $i$ denotes the inclusion map $i \colon P_{\le k} \hookrightarrow P$.
\end{proposition}

\begin{proof}
Write $Q = P_{\le k}$, $\Gamma = \Lv$, $H = F|_{P_{\le k}}$ and $K = G|_{P_{\le k}}$. We have the following commutative diagram
\begin{center}
\begin{tikzpicture}
    \diagram{d}{3em}{3em}{
        Q & P \\
        Q & P \\
    };

    \path[->,font = \scriptsize, midway]
    (d-1-1) edge node[midway, right]{$\Gamma|_Q$} (d-2-1)
    (d-1-1) edge node[midway, above]{$i$} (d-1-2)
    (d-2-1) edge node[midway, above]{$i$} (d-2-2)
    (d-1-2) edge node[midway, right]{$\Gamma$} (d-2-2);
\end{tikzpicture}
\end{center}
We want to show that given a $\Gamma|_Q$--interleaving $(\phi, \psi)$ between $H$ and $K$, there is a $\Gamma$--interleaving $(\hat \phi, \hat \psi)$ between the left Kan extensions of $H$ and $K$ along the inclusion $i \colon Q \hookrightarrow P$.

Let $\Lke H$ be the left Kan extension of $H$ along $i$, and let $\Lke K$ be the left Kan extension of $K$ along $i$. Let $\alpha \colon H \to \Lke H \circ i$ and $\beta \colon K \to \Lke K \circ i$ be the corresponding natural transformations.

We define $\hat \phi \colon \Lke H \to \Lke K \circ \Gamma$ by applying the universal property of $\Lke H$ on the following diagram.
\begin{center}
\begin{tikzpicture}
    \diagram{d}{3em}{3em}{
        \ & \Lke H \circ i & \ \\
        H & K \circ \Gamma|_Q & \Lke K \circ i \circ \Gamma|_Q = \Lke K \circ \Gamma \circ i\\
    };

    \path[->,font = \scriptsize, midway]
    (d-2-1) edge node[midway, left]{$\alpha$} (d-1-2)
    (d-2-2) edge node[midway, above]{$\beta_{\Gamma|_Q}$} (d-2-3)
    (d-2-1) edge node[midway, above]{$\phi$} (d-2-2);
    \path[->,font = \scriptsize, midway, dashed]
    (d-1-2) edge node[midway, above]{$\hat \phi_i$} (d-2-3);
\end{tikzpicture}
\end{center}
We define $\hat \psi \colon \Lke K \to \Lke H \circ \Gamma$ similarly.

We now need to show that ${\hat \psi}_{\Gamma} \circ \hat \phi = \Lke H(\eta_{(\Gamma|_Q)^2})$ and ${\hat \phi}_{\Gamma} \circ \hat \psi = \Lke K(\eta_{(\Gamma|_Q)^2})$.
We prove the first equality. The proof of the second equality is similar. Consider the diagram
\begin{center}
\begin{tikzpicture}
    \diagram{d}{3em}{3em}{
        \Lke H \circ i & \Lke K \circ \Gamma \circ i = \Lke K \circ i \circ \Gamma|_Q  & \Lke H \circ \Gamma \circ i \circ \Gamma|_Q = \Lke H \circ i \circ (\Gamma|_Q)^2 \\
        H & K \circ \Gamma|_Q & H \circ (\Gamma|_Q)^2\\
    };

    \path[->,font = \scriptsize, midway]
    (d-2-1) edge node[midway, left]{$\alpha$} (d-1-1)
    (d-2-2) edge node[midway, left]{$\beta_{\Gamma|_Q}$} (d-1-2)
    (d-2-3) edge node[midway, left]{$\alpha_{(\Gamma|_Q)^2}$} (d-1-3)
    (d-2-1) edge node[midway, above]{$\phi$} (d-2-2)
    (d-2-2) edge node[midway, above]{$\psi_{\Gamma|_Q}$} (d-2-3)
    (d-2-1) edge node[midway, above]{$\phi$} (d-2-2)
    (d-2-2) edge node[midway, above]{$\psi_{\Gamma|_Q}$} (d-2-3)
    (d-1-1) edge node[midway, above]{${\hat \phi}_i$} (d-1-2)
    (d-1-2) edge node[midway, above]{${\hat \psi}_{i \Gamma |_Q}$} (d-1-3);
\end{tikzpicture}
\end{center}
As both squares commute, the entire diagram commutes. The top-left composite is
\begin{equation*}
    {\hat \psi}_{i \Gamma |_Q} \circ {\hat \phi}_i \circ \alpha = {\hat \psi}_{\Gamma i} \circ {\hat \phi}_i \circ \alpha = ({\hat \psi}_{\Gamma} \circ \hat \phi)_i \circ \alpha,
\end{equation*}
while the bottom-right composite is
\begin{equation*}
    \alpha_{(\Gamma|_Q)^2} \circ \psi_{\Gamma|_Q} \circ \phi
    = \alpha_{(\Gamma|_Q)^2} \circ H(\eta_{(\Gamma|_Q)^2}).
\end{equation*}
We further have that $\Lke H(\eta_{\Gamma^2})_i \circ \alpha = \alpha_{(\Gamma|_Q)^2} \circ H(\eta_{(\Gamma|_Q)^2})$. To see this, evaluate both sides in $x \in Q$:
\begin{align*}
    \left( \Lke H(\eta_{\Gamma^2})_i \circ \alpha \right)_x &=
    \Lke H \left( i(x) \le \Gamma^2(i(x)) \right) \circ \alpha_x \\
    &= \Lke H\left( i(x) \le i((\Gamma|_Q)^2(x)) \right) \circ \alpha_x \\
    &= (\Lke H \circ i)\left( x \le (\Gamma|_Q)^2(x) \right) \circ \alpha_x \\
    &= \alpha_{(\Gamma|_Q)^2 (x)} \circ H \left(x \le (\Gamma|_Q)^2(x) \right) \qquad \text{(by the naturality of $\alpha$)} \\
    &= \left(\alpha_{(\Gamma|_Q)^2} \circ H(\eta_{(\Gamma|_Q)^2})\right)_x.
\end{align*}
In conclusion, we get that 
$({\hat \psi}_{\Gamma} \circ \hat \phi)_i \circ \alpha = \alpha_{(\Gamma|_Q)^2} \circ H(\eta_{(\Gamma|_Q)^2}) = \Lke H(\eta_{\Gamma^2})_i \circ \alpha$. By the universal property of left Kan extensions, this gives that ${\hat \psi}_{\Gamma} \circ \hat \phi = \Lke H(\eta_{\Gamma^2})$, as desired.

\end{proof}

\begin{lemma}\label{lemma:interleaving_cof_repl}
    Let $F, G \colon \R_{\ge 0}^n \to \M$ be functors to a good model category, and let $Z \colon \Ib^{\Lv} \to \M$ be a $\Lv$--interleaving from $F$ to $G$. Let $QZ$ be a cofibrant replacement of $Z$. Then $QZ$ is an interleaving between cofibrant replacements of $F$ and $G$.
\end{lemma}

\begin{proof}
It's clear that $(QZ) \circ E^0 \simeq F$ and $(QZ) \circ E^1 \simeq G$. Thus, it suffices to show that $(QZ) \circ E^0$ and  $(QZ) \circ E^1$ are cofibrant.
We show that $(QZ) \circ E^0$ is cofibrant. The proof for $(QZ) \circ E^1$ is similar.

Let $r \colon \Ib^{\Lv} \to \R_{\ge 0}^n$ be the map given by
\begin{equation*}
    r(x,i) =
    \begin{cases}
        x, &\quad i=0, \\ 
        \Lv(x), &\quad i=1. \\ 
    \end{cases}
\end{equation*}
By checking the four cases for values of $i$, it is easily verified that this map is order-preserving, and hence a functor.

We claim that we have an adjunction $r \dashv E^0$. We show this by finding a counit and unit (commutativity conditions are trivially satisfied as natural transformations between poset-valued functors are unique if they exist). As $r \circ E^0$ is the identity, it suffices to find a unit $\id \to E^0 \circ r$. In other words, we need to show that
\begin{equation*}
    (E^0 \circ r)(x, i) \ge (x,i)
\end{equation*}
for all $(x,i) \in \Ib^{\Lv}$. In the case of $i=0$, we get
\begin{equation*}
    (E^0 \circ r)(x, 0) = E^0(x) = (x,0),
\end{equation*}
so the relation is trivially satisfied. If $i=1$, we get
\begin{equation*}
    (E^0 \circ r)(x, 1) = E^0(\Lv(x)) = (\Lv(x), 0) \ge (x, 1).
\end{equation*}
Thus, $r \dashv E^0$, as desired

Now, it follows that
\[
    \begin{tikzcd}[column sep=large]
        \Fun(\Ib^{\Lv}, \C) \arrow[r, shift left=1ex, "(E^0)^*"{name=G, yshift=1pt}] & \Fun(\R_{\ge 0}^n, \C) \arrow[l, shift left=.5ex, "r^*"{name=F}]
        \arrow[phantom, from=F, to=G, , "\scriptscriptstyle\boldsymbol{\bot}"],
    \end{tikzcd}
\]
is an adjunction. Furthermore, as $r^*$ preserves fibrations and weak equivalences, the adjunction is a Quillen adjunction. Hence, $(E^0)^*$ preserves cofibrant objects, and we are done.
\end{proof}

\begin{corollary}\label{cor:stability_theorem}
    Let $F, G \colon \R_{\ge 0}^n \to \M$ be functors to a good model category, and suppose that $F$ and $G$ are $\Lv$--interleaved. Then $T_k F$ and $T_k G$ are $\Lv$--homotopy-interleaved.
\end{corollary}

\begin{proof}
Let $P = \R_{\ge 0}^n$, and
let $i$ denote the inclusion map $i \colon P_{\le k} \hookrightarrow P$.
By \autoref{lemma:interleaving_cof_repl}, we have a $\Lv$--interleaving between $QF$ and $QG$, where $QF$ and $QG$ are cofibrant replacements of $F$ and $G$, respectively.
Furthermore, by \autoref{prop:interleaving_lke}, we have a $\Lv$--interleaving between $\Lan_i (QF|_{P_{\le k}})$ and $\Lan_i (QG|_{P_{\le k}})$.
Finally, we have $T_k F \simeq \Lan_i (QF|_{P_{\le k}})$ and $T_k G \simeq \Lan_i (QG|_{P_{\le k}})$, as desired.
\end{proof}

\begin{corollary}
    Let $F, G \colon (\R_{\ge 0})^n \to \M$ be functors to a good model category. For all $k \in \Z_+$, we have
    \begin{equation*}
        d^{\Lambda}_H (T_k F, T_k G) \le d^{\Lambda}_H (F, G)
    \end{equation*}
    and
    \begin{equation*}
        d^{\Lambda}_H (T_k F, T_k G) \le d^{\Lambda}_H (T_{k+1} F, T_{k+1} G).
    \end{equation*}
\end{corollary}

\begin{proof}
The first inequality follows directly from \autoref{cor:stability_theorem}.
The second follows from \autoref{cor:stability_theorem} and the fact that $T_k F \simeq T_k (T_{k+1} F)$.
\end{proof}

The preceding corollary can be summarized in the following sequence of inequalities.

\begin{equation*}
    d^{\Lambda}_H (T_1 F, T_1 G) \le d^{\Lambda}_H (T_2 F, T_2 G) \le \dots \le d^{\Lambda}_H (F, G)
\end{equation*}

\section{Example: functors between posets}

We consider here functors of the form $F \colon P \to Q$, where $P$ and $Q$ are posets.

\subsection{Notation}

For a set $V$, we let $\Po(V)$ denote the power set of $V$. For a poset $P$, we let $P^{\op}$ denote the opposite poset.

For a simplicial complex $X$, we let $X_n$ denote its set of $n$--simplices. In particular, $X_0$ is the vertex set of $X$. The \emph{$n$--skeleton} of $X$ is defined as 
\begin{equation*}
    X_{\le n} = \bigcup_{i = 0}^n X_i = \{\sigma \in X : |\sigma| \le n+1\}.
\end{equation*}
Given a subcomplex $Y \subset X$, we write $X \searrow Y$ to denote that $Y$ can be reached from $X$ through a sequence of simplicial collapses.

A simplicial complex $X$ is called \emph{$n$--skeletal} if it has no $i$--simplices for $i > n$. A simplicial complex $X$ is called \emph{$n$--coskeletal} if, for each $i > n$, it has a single $i$--simplex for each compatible set of faces. In other words, $X$ is $n$--coskeletal if it is the maximal simplicial complex having $X_{\le n}$ as its $n$--skeleton (maximal in the sense of having the most simplices).

\subsection{The poset of simplicial complexes over a vertex set}

Recall that $\Po(V)$ is a join-factorization lattice when $V$ is finite: It is distributive because $A \cap (B \cup C) = (A \cap B) \cup (A \cap C)$, and the indecomposable join-decomposition of a finite set $\{x_1, \dots, x_k\}$ is $\{\{x_1\}, \dots, \{x_k\}\}$. The join-dimension of an element $U \in \Po(V)$ is the cardinality of $U$. Hence, an $n$--simplex has join-dimension $n+1$.

Let $V$ be a set of vertices. A \emph{simplicial complex over $V$} is a subset $X \subseteq \Po(V)$ subject to the following two conditions:
\begin{itemize}
    \item $\emptyset \in X$.
    \item If $\sigma \subseteq \tau$ and $\tau \in X$, then $\sigma \in X$.
\end{itemize}

We can give a second, equivalent formulation as follows.
Equip $\Po(V)$ with the partial order $\sigma \le \tau \Leftrightarrow \sigma \subseteq \tau$, and let $\{0, 1\}$ be the poset with two objects and the partial order $0 < 1$.
A \emph{simplicial complex over $V$} is a functor
\begin{equation}\label{eq:def_simpl_cplx}
    X \colon \Po(V)^{\op} \to \{0, 1\},
\end{equation}
that satisfies $X(\emptyset) = 1$. 

To see that this second formulation is equivalent to the first one, note that a subset $X \subseteq \Po(V)$ is canonically identified with a function $f \colon \Po(V) \to \{0, 1\}$ (where $f(v) = 1 \Leftrightarrow v \in X$). Furthermore, the functoriality condition in \eqref{eq:def_simpl_cplx} says precisely that $\tau \supseteq \sigma \Rightarrow X(\tau) \le X(\sigma)$, or in other words, if $\tau \supseteq \sigma$ and $X(\tau) = 1$, then $X(\sigma) = 1$. Finally, the condition $X(\emptyset) = 1$ simply corresponds to $\emptyset \in X$.

We now define the \emph{poset of simplicial complexes over $V$} as the subcategory of the functor category
\begin{equation}
    \SCplx_V \subset \Fun(\Po(V)^{\op}, \{0,1\})
\end{equation}
consisting of functors sending $\emptyset$ to 1.
The functor category $\Fun(\Po(V)^{\op}, \{0,1\})$ is a poset as the category of functors between two posets is a poset. Explicitly, given two simplicial complexes $X, Y \colon \Po(V)^{\op} \to \{0,1\}$, we have $X \le Y$ if and only if $X(\sigma) \le Y(\sigma)$ for all $\sigma \in \Po(V)$. Using the first formulation, this is the same as saying that $X$ is a subset (or \emph{subcomplex}) of $Y$.

The poset $\SCplx_V$ is also a lattice. As $\{0, 1\}$ is a finite lattice, and hence bicomplete, then so is $\Fun(\Po(V)^{\op}, \{0,1\})$.  Given two simplicial complexes $X, Y \colon \Po(V)^{\op} \to \{0,1\}$, their meet is given by $(X \wedge Y)(\sigma) = X(\sigma) \wedge Y(\sigma)$, and their join is given by $(X \vee Y)(\sigma) = X(\sigma) \vee Y(\sigma)$. Under the first formulation, the meet and join correspond to intersection and union of sets, respectively. In other words, $X \vee Y = X \cup Y$ and $X \wedge Y = X \cap Y$.
Finally, it is easily that the subposet $\SCplx_V$ is closed under meet and join.

\begin{example}[Codegree captures the concept of $n$--skeletal]\label{ex:n_skeletal}
Let $X \in \SCplx_V$ be a simplicial complex over $V$.
Consider the functor
\begin{equation*}
    F_X \colon \Po(V) \to \SCplx_V
\end{equation*}
that sends a subset $U \subseteq V$ to the subcomplex of $X$ spanned by the vertices $U$. If $U \subseteq U' \subseteq V$, then clearly $F_X(U)$ is a subcomplex of $F_X(U')$, and hence $F_X$ is order-preserving (and thus a well-defined functor).

Equivalently, we can write
\begin{equation*}
    F_X(U) = \Po(U) \wedge X = \Po(U) \cap X,
\end{equation*}
where we consider $\Po(U) \subseteq \Po(V)$ as a simplicial complex.

I will now show that $F_X$ is codegree $n$ if and only if $X$ is $(n-1)$--skeletal. Using \autoref{theoremA_cat} and \autoref{theoremB_cat}, $F_X$ is codegree $n$ if and only if $F_X = T_n F_X$. Given a set $U \subseteq V$,
\begin{align*}
    T_n F_X(U) &= \underset{\sigma \subseteq U, |\sigma| \le n}{\colim} F_X(\sigma) \\
    &= \bigcup_{\sigma \subseteq U, |\sigma| \le n} (\Po(\sigma) \cap X) \\
    &= X \cap \bigcup_{\sigma \subseteq U, |\sigma| \le n} \Po(\sigma) \\
    &= X \cap \left(\Po(U)\right)_{\le n-1} = \left(X \cap \Po(U)\right)_{\le n-1}.
\end{align*}
In words, $T_n F_X (U)$ is the subset of $X_{\le n-1}$ spanned by the vertices in $U$. Thus, if $X$ is $(n-1)$--skeletal, then $T_n F_X (U) = F_X (U)$ for all $U \subseteq V$. Conversely, if $T_n F_X = F_X$, then $X = F_X(V) = T_n F_X (V) = X_{\le n-1}$.
\end{example}

\begin{remark}
In \autoref{ex:n_skeletal} we used \autoref{theoremA_cat} and \autoref{theoremB_cat} to show that $F_X$ is codegree $n$ when $X$ is $(n-1)$--skeletal.
We now get ``for free'' that when $X$ is $(n-1)$--skeletal, then $F_X$ sends strongly bicartesian $(n+1)$--cubes to cocartesian $(n+1)$--cubes. We explore the consequences of this.

For $n=1$, i.e., when $X$ is 0--skeletal, we simply get that
\begin{equation*}
    F_X(A \cup B) = F_X(A) \cup F_X(B).
\end{equation*}
This makes sense, as $X$ consists only of 0--simplices and $F_X(U) = \{\{v\} : v \in U\}$.

For $n=2$, $X$ is 1--skeletal, so a graph.
The codegree $n$ condition now tells us that when $U_1, U_2, U_3$ is a pairwise cover of $V$ (i.e., $U_1 \cup U_2 = U_1 \cup U_3 = U_2 \cup U_3 = V$), then
\begin{equation*}
    X = F_X(V) = F_X(U_1) \cup F_X(U_2) \cup F_X(U_3).
\end{equation*}
This equation essentially tells us that if you know the subgraphs for three such vertex sets, then you can reconstruct the entire graph. Knowing only two such subgraphs would not be enough! Consider the graph in figure \autoref{fig:triangle}. Knowing the subgraphs spanned by $\{a, b\}$ and $\{a, c\}$ would not be enough to determine the entire graph, but if you know the subgraph spanned by $\{b, c\}$ as well, you have enough information.

\begin{figure}[htbp]
    \centering
    \begin{tikzpicture}[scale=.8]
        % Edges
        \draw (4, 0) -- (2, 3.4);
        \draw (4, 0) -- (6, 3.4);
        \draw (2, 3.4) -- (6, 3.4);
        
        % Vertices
        \filldraw (2, 3.4) circle (2 pt) node[anchor=east]{$a$};
        \filldraw (4, 0) circle (2 pt) node[anchor=east]{$b$};
        \filldraw (6, 3.4) circle (2 pt) node[anchor=west]{$c$};
    \end{tikzpicture}
    \caption{A graph.}
    \label{fig:triangle}
\end{figure}

Note that finding a pairwise cover of $V$ is easy; one can partition $V$ into three sets $A,B,C$ and take $U_1 = A \cup B, U_2 = A \cup C$ and $U_3 = B \cup C$. Consequently, if $A, B, C$ are three sets so that $A \cup B \cup C = V$, then
\begin{equation*}
    X = F_X(V) = F_X(A \cup B) \cup F_X(A \cup C) \cup F_X(B \cup C).
\end{equation*}

The same example works for higher $n$, in which case one would need to take a pairwise cover consisting of $n+1$ sets.
\end{remark}

\begin{example}[Codegree captures the concept of $n$--coskeletal]\label{ex:n_coskeletal}
Given a simplicial complex $X \colon \Po(V)^{\op} \to \{0, 1\}$, consider its opposite functor $X^{\op} \colon \Po(V) \to \{0,1\}^{\op}$. I show that $X^{\op}$ is codegree $n$ if and only if $X$ is $(n-1)$--coskeletal.

Using \autoref{theoremA_cat} and \autoref{theoremB_cat}, $X^{\op}$ is codegree $n$ precisely if for all $\sigma \subseteq V$ with $|\sigma| > n$,
\begin{equation*}
    X^{\op}(\sigma) = \underset{\tau \subseteq \sigma, |\tau| \le n}{\colim} X^{\op} (\tau) = \bigvee_{\tau \subseteq \sigma, |\tau| \le n} X^{\op} (\tau)
\end{equation*}
This is equivalent to
\begin{equation*}
    X(\sigma) = \bigwedge_{\tau \subseteq \sigma, |\tau| \le n} X (\tau),
\end{equation*}
or, in other words, that $X(\sigma) = 1$ if and only if $X(\tau) = 1$ for all $\tau \le \sigma$ with $|\tau| \le n$. This is precisely saying that $X$ is coskeletal.
\end{example}

\begin{remark}
As for \autoref{ex:n_skeletal}, we can use \autoref{def:deg_n_bicartesian_cat} together with \autoref{ex:n_coskeletal} to get a property that is satisfied by $(n-1)$--coskeletal simplicial complexes.

Let $X \in \SCplx_V$ be $(n-1)$--coskeletal and $\tau_1, \dots, \tau_{n+1}$ a size $n$ pairwise cover of a set $\sigma \subseteq V$. Unwinding \autoref{def:deg_n_bicartesian_cat}, we get
\begin{equation}\label{eq:coskel_prop}
    X(\sigma) = \bigwedge_{i=1}^{n+1} X(\tau_i).
\end{equation}

In the case of $n=1$, this simply says that for a 0--coskeletal simplicial complex $X$, if $\tau_1 \in X$ and $\tau_2 \in X$, then $(\tau_1 \cup \tau_2) \in X$.

In the $n=2$ case, \eqref{eq:coskel_prop} means that the following holds when $X$ is 1--coskeletal. Let $\sigma \subseteq V$ and $\tau_1, \tau_2, \tau_3 \subseteq \sigma$ such that $\tau_1 \cup \tau_2 = \tau_1 \cup \tau_3 = \tau_2 \cup \tau_3 = \sigma$. If $\tau_1, \tau_2$ and $\tau_3$ are all in $X$, then $\sigma$ is in $X$. This is illustrated in \autoref{fig:tetrahedron_coskeletal}.

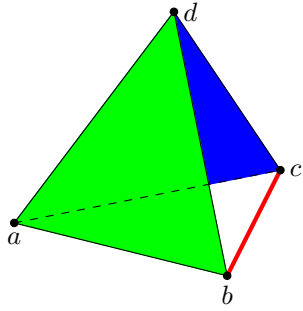
\begin{figure}[htbp]
    \centering
    \begin{tikzpicture}[scale=.7]
        % Triangles
        \filldraw[blue,opacity=.5] (0,1) -- (5, 2) -- (3, 5) -- cycle;
        \filldraw[green,opacity=.5] (0,1) -- (4,0) -- (3, 5) -- cycle;
        
        % Edges
        \draw (0, 1) -- (4, 0);
        \draw (0, 1) -- (3, 5);
        \draw (4, 0) -- (3, 5);
        \draw[red, ultra thick] (4, 0) -- (5, 2);
        \draw (3, 5) -- (5, 2);
        \draw[dashed] (0, 1) -- (5, 2);
        
        % Vertices
        \filldraw (0, 1) circle (2 pt) node[anchor=north]{$a$};
        \filldraw (4, 0) circle (2 pt) node[anchor=north]{$b$};
        \filldraw (5, 2) circle (2 pt) node[anchor=west]{$c$};
        \filldraw (3, 5) circle (2 pt) node[anchor=west]{$d$};
    \end{tikzpicture}
    \caption{Illustrated is a 3--simplex $\sigma$. Suppose that $X$ is a 1--coskeletal simplicial complex and $\sigma \in X$. If $\{a,b,d\}, \{a,c,d\}$ and \{b,c\} are all in $X$, then $\{a,b,c,d\}$ is in $X$.}
    \label{fig:tetrahedron_coskeletal}
\end{figure}
\end{remark}

\begin{example}[Codegree of a filtration]\label{ex:deg_filtration}
Let $V$ be a vertex set.
Given a functor $f \colon \Po(V) \to \nnR$, we can define a \emph{filtration} on $\Po(V)$ as
\begin{align*}
    F_f \colon \nnR &\to \SCplx_V \\
    t &\mapsto \{\sigma : f(\sigma) \le t\}.
\end{align*}

For an expression $P$, let $[P]$ denote its truth value in $\{0,1\}$.
Suppose that $f \colon \Po(V) \to \nnR$ is a codegree $n$ functor.\footnote{Technically, $\nnR$ is not a cocomplete category as it doesn't have a maximal element, and therefore we cannot apply poset cocalculus in this situation. However, as the source poset is finite here, this doesn't cause any issues. We could choose to work with the bicomplete poset $[0, \infty]$ instead, and all the computations in this example would be the same.}
Then, for $t \in \nnR$ and $\sigma \in \Po(V)$,
\begin{align*}
    F_f (t) (\sigma) &= \left[f(\sigma) \le t\right] \\
    &= \left[\left(\bigvee_{\tau \subseteq \sigma, |\tau| \le n} f(\tau)\right) \le t\right] \\
    &= \bigwedge_{\tau \subseteq \sigma, |\tau| \le n} \left[f(\tau) \le t\right] \\
    &= \bigwedge_{\tau \subseteq \sigma, |\tau| \le n} F_f(t)(\tau).
\end{align*}
Hence, when $f$ is codegree $n$, then the filtration $F_f(t)$ is $(n-1)$--coskeletal for all $t \in \nnR$.
An example of a filtration satisfying this property is the \emph{Vietoris-Rips filtration}. Indeed, each complex in a Vietoris-Rips filtration is 1--coskeletal (commonly known as a \emph{clique complex}).

The implication goes both ways. If $F_f(t)$ is $(n-1)$--coskeletal at every $t$, then for every $t \in \R$ and $\sigma \in \Po(V)$,
\begin{align*}
    F_f (t) (\sigma) &= \bigwedge_{\tau \subseteq \sigma, |\tau| \le n} F_f(t)(\tau) \\
    &= \bigwedge_{\tau \subseteq \sigma, |\tau| \le n} F_{T_n f}(t)(\tau) \\
    &= F_{T_n f}(t)(\sigma).
\end{align*}
Hence, $F_f = F_{T_n f}$, and so $f = T_n f$.
\end{example}

Let $V \in \R^n$ be a finite point cloud.
The \emph{Čech filtration} of $V$ is the filtration $F_{f_\mathbf{C}}$ induced by the functor
\begin{align*}
    f_{\mathbf{C}} \colon \Po(V) &\to \R_{\ge 0} \\
    U &\mapsto \inf\left\{\varepsilon \ : \ \exists z \in \R^n \textrm{ such that } U \subseteq \overline{B_{\varepsilon}}(z) \right\}.
\end{align*}
In other words, $f_{\mathbf{C}}(U)$ is the radius of $U$.

The \emph{Vietoris-Rips filtration} of $V$ is the filtration $F_{f_\mathbf{VR}}$ induced by the functor
\begin{align*}
    f_{\mathbf{VR}} \colon \Po(V) &\to \R_{\ge 0} \\
    U &\mapsto \max\left\{ d(x, y) : x,y \in U \right\}.
\end{align*}
In other words, $f_{\mathbf{VR}}(U)$ is the diameter of $U$.

The following proposition shows that the Vietoris-Rips filtration is in some sense a codegree 2 approximation of the Čech filtration, scaled by a factor.

\begin{proposition}\label{prop:vietoris_rips_cech}
    Let $V \subset \R^d$ be a finite point cloud, and let $f_{\mathbf{C}} \colon \Po(V) \to \nnR$ and $f_{\mathbf{VR}} \colon \Po(V) \to \nnR$ be the Čech filtration function and the Vietoris-Rips filtration function, respectively.
    Then,
    \begin{equation*}
        f_{\mathbf{VR}} = 2 \cdot T_2 f_{\mathbf{C}}.
    \end{equation*}
\end{proposition}

\begin{proof}
It is straightforward to check that $f_{\mathbf{VR}}$ is codegree 2:
\begin{align*}
    \left(T_2 f_{\mathbf{VR}}\right)(U) &= \bigvee_{W \subseteq U, |W| \le 2} f_{\mathbf{VR}}(W) \\
    & = \bigvee_{\{x,y\} \subseteq U} f_{\mathbf{VR}}(x,y) \\
    & = \bigvee_{\{x,y\} \subseteq U} d(x,y) = f_{\mathbf{VR}}(U). \\
\end{align*}

Furthermore, we have that $f_{\mathbf{VR}}(U) = 2 \cdot f_{\mathbf{C}}(U)$ when $|U| \le 2$. When $U$ is the empty set or a singleton, both functors simply evaluate to 0. When $U = \{x,y\}$, $f_{\mathbf{C}}(\{x,y\}) = d(x,y)/2$, as desired. Now, as $f_{\mathbf{VR}}$ is codegree $2$, we get from \autoref{theoremB_cat} that
\begin{equation*}
    f_{\mathbf{VR}} = T_2 \left(f_{\mathbf{VR}}\right) = T_2 \left(2 \cdot f_{\mathbf{C}}\right) = 2 \cdot T_2 f_{\mathbf{C}},
\end{equation*}
(where the last equality follows from the fact that multiplication distributes over supremum in $\R_{\ge 0}$).
\end{proof}

\begin{example}
Let $V \subset \R^d$ be a finite set of points, and let $f_{\mathbf{C}} \colon \Po(V) \to \nnR$ be its corresponding Čech filtration. It is a known consequence of Helly's theorem that $F_{f_{\mathbf{C}}}(t)$ is $d$--coskeletal at every $t \in \nnR$.
% Possible TODO: maybe add source here?
Hence, $f_{\mathbf{C}}$ is codegree $d+1$.
\end{example}

\begin{remark}
Let $\iota$ denote the functor
\begin{align*}
    \iota \colon \Po(V) &\to \SCplx_V \\
    U &\mapsto \Po(U).
\end{align*}
Let $f \colon \Po(V) \to \R$. The filtration $F_f$ in \autoref{ex:deg_filtration} is the left Kan extension of $\iota$ along $f$.
\end{remark}

\section{Dual calculus}

The constructions in poset cocalculus can be dualized to create a poset \emph{calculus}, which we sketch here. Let $P$ be a join-factorization lattice and let $\M$ be a model category such that $\Fun(P^{\op}, \M)$ admits the injective model structure. The latter condition is equivalent to asking that $\Fun(P, \M^{\op})$ admits the projective model structure. Hence, for a functor $F \colon P^{\op} \to \M$ we can construct the Taylor telescope for $F^{\op} \colon P \to \M^{\op}$. Taking the opposite of this telescope gives a tower of functors $P^{\op} \to \M$,
\begin{equation*}
\begin{aligned}
\begin{tikzpicture}
    \diagram{d}{3em}{3em}{
        \ & \vdots \\
        \ & T^2 F \\
        F & T^1 F, \\
    };
    
    \path[->,font = \scriptsize, midway]
    (d-1-2) edge node[left]{} (d-2-2)
    (d-3-1) edge (d-1-2)
    (d-2-2) edge node[left]{} (d-3-2)
    (d-3-1) edge node[below]{} (d-2-2)
    (d-3-1) edge node[below]{} (d-3-2);
\end{tikzpicture}
\end{aligned}
\end{equation*}
where we denote by $T^n F$ the $n$th functor in this tower (i.e., we use superscript instead of subscript for the dual case). We call $T^n F$ the \emph{degree $n$ approximation} of $F$. As a homotopy colimit in $\M^{\op}$ is a homotopy limit in $\M$, it follows that $T^n F$ sends strongly bicartesian $(n+1)$--cubes to homotopy \emph{cartesian} $(n+1)$--cubes. We similarly get dual versions of \autoref{theoremA_cat} and \autoref{theoremB_cat}. The formula for $T^n F$ is
\begin{equation*}
    T^n F(x) = \underset{v \in P_{\le n}, v \le x}{\holim} F(v),
\end{equation*}
where we let $\le$ be the partial order relation in $P$ (note: it's \emph{not} the partial order relation in $P^{\op}$).

\subsection{Examples}

\subsubsection{$n$--coskeletal simplicial complexes}

\begin{example}
We revisit the example in \autoref{ex:n_coskeletal}. Recall that we had to consider the opposite of the functor $X \colon \Po(V)^{\op} \to \{0,1\}$. With the dual framework, poset calculus, we can consider $X$ directly.

Let $X \colon \Po(V)^{\op} \to \{0,1\}$ be a simplicial complex. Then
\begin{align*}
    T^n X (\sigma) &= \lim_{\tau \le \sigma, |\tau| \le n} X(\tau) \\
    &= \bigwedge_{\tau \le \sigma, |\tau| \le n} X(\tau).
\end{align*}
As in \autoref{ex:n_coskeletal}, we see that $X$ is $(n-1)$--coskeletal precisely if $X = T^n X$, which holds precisely when $X$ sends strongly bicartesian $(n+1)$--cubes to cartesian $(n+1)$--cubes.

We can also show this directly from the definition of $T^n X$. By definition, $T^n X = (T_n X^{\op})^{\op}$, so $T^n X = X$ if and only if $T_n X^{\op} = X^{\op}$.
\end{example}

\subsubsection{Internal hom of simplicial complexes}

\begin{example}
We define a \emph{simplicial map} between two simplicial complexes $X$ and $Y$ as a function $f \colon X_0 \to Y_0$ between the vertex sets such that $f(\sigma) \in Y$ for every simplex $\sigma \in X$.
Let $\SCplx$ be the category of simplicial complexes and simplicial maps. This category is cartesian closed (\cite[1.4]{scplxshom}, \cite[IV.7]{closedCategories}),
where the tensor product is the categorical product. Explicitly, for $X, Y \in \SCplx$,
\begin{itemize}
    \item the vertex set of $X \times Y$ is $(X \times Y)_0 = X_0 \times Y_0$, and
    \item $U = \{(v_0, w_0), \dots, (v_n, w_n)\}$ is in $X \times Y$ if and only if $\pi_X(U) \in X$ and $\pi_Y(U) \in Y$, where $\pi_X \colon X\times Y \to X$ and $\pi_Y \colon X \times Y \to Y$ are the projections.
\end{itemize}
The internal hom functor
\begin{equation*}
    \Hhom \colon \SCplx^{\op} \times \SCplx \to \SCplx,
\end{equation*}
is defined as follows.
\begin{itemize}
    \item The vertices in $\Hhom(X,Y)$ are the simplicial maps from $X$ to $Y$, i.e.,
    \begin{equation*}
        \Hhom(X,Y)_0 = \Hom(X,Y).
    \end{equation*}
    
    \item The $n$--simplex $\{f_0, \dots, f_n\}$ is in $\Hhom(X,Y)$ if and only if
    \begin{equation*}
        \bigcup_{i=0}^n f_i(\sigma) \in Y \ \textrm{ for every } \sigma \in X.
    \end{equation*}
\end{itemize}
The Hom-tensor adjunction is defined on 0--simplices as
\begin{align*}
    \Hhom(X \times Y, Z)_0 &\to \Hhom(X, \Hhom(Y, Z))_0 \\
    f &\mapsto g
\end{align*}
where $g(x)(y) = f(x,y)$. For more details on the internal hom functor for simplicial complexes, see \cite[1.4]{scplxshom} and \cite[IV.7]{closedCategories}.
% maybe TODO: define realization of simplicial complexes (=|N(P \ {Ø})|)

Suppose we are given two simplicial complexes $X,Y \in \SCplx$. We want to study the homotopy type of $|\Hhom(X,Y)|$, where $| - | \colon \SCplx \to \Spaces$ is the functor that takes a simplicial complex to its geometric realization, and $\Spaces$ is the category of compactly generated topological spaces.

Let $V = X_0$ be the vertex set of $Y$. Recall that $\SCplx_V$ denotes the poset of subcomplexes of $\Po(V)$. Consider the functor
\begin{align*}
    F \colon \SCplx_V &\to \Spaces \\
    Z &\mapsto | \Hhom(Z, Y) |.
\end{align*}

We now show that $F$ satisfies
\begin{equation}\label{eq:hom_example}
    T^1 F (Z) \simeq \Hhom(|Z|, |Y|),
\end{equation}
where $\Hhom$ denotes the internal hom functor in $\Spaces$.

First, we denote by $G$ the functor
\begin{align*}
    G \colon \SCplx_V^{\op} &\to \Spaces \\
    Z &\mapsto \Hhom(|Z|, |Y|).
\end{align*}
We can visualize $F$ and $G$ as the top-right and bottom-left compositions, respectively, in the following diagram.
\begin{center}
\begin{tikzcd}
\SCplx_V^{\op} \arrow[r, "i^{\op}", hook] & \SCplx^{\op} \arrow[rr, "{\Hhom(-, Y)}"] \arrow[d, "| - |^{\op}"'] & & \SCplx \arrow[d, "| - |"] \\
    & \Spaces^{\op} \arrow[rr, "{\Hhom(-, | Y |)}"'] & & \Spaces                  
\end{tikzcd} 
\end{center}
Here, $i$ is the inclusion of $\SCplx_V$ into $\SCplx$.

We will show that $G$ sends bicartesian squares to homotopy cartesian squares (\textbf{I}) and that $G(Z) \simeq F(Z)$ when $\jdim(Z) \le 1$. This will imply \eqref{eq:hom_example}, by  the dual of \autoref{theoremB}.

\begin{description}[style=multiline]
    \item[I] We show that $G$ sends strongly bicartesian squares to homotopy cartesian squares. A strongly bicartesian square in $\SCplx_V$ is of the form
    \begin{center}
    \begin{tikzcd}
        Z \cap W \arrow[r] \arrow[d] & Z \arrow[d] \\
        W \arrow[r]                  & Z \cup W   
    \end{tikzcd}
    \end{center}
    One can verify directly that
    \begin{center}
    \begin{tikzcd}
        \mid Z \cap W\mid \arrow[r] \arrow[d] & \mid Z \mid \arrow[d] \\
        \mid W \mid  \arrow[r] &  \mid Z \cup W \mid
    \end{tikzcd}
    \end{center}
    is cocartesian. Furthermore, the top and left arrows are cofibrations, so the square is homotopy cocartesian.

    Furthermore, $\Spaces$ is a symmetric monoidal model category \cite[Proposition 4.2.11]{hovey}, so $\Hhom(- , |Y|)$ is right Quillen (as a functor $\Spaces^{\op} \to \Spaces$).
    Therefore, $\Hhom(- , |Y|)$ sends colimits to limits and cofibrations to fibrations.
    In conclusion, the square
    \begin{center}
    \begin{tikzcd}
        \Hom(\mid Z \cap W \mid, \mid Y \mid)  & \Hom(\mid Z \mid, \mid Y \mid) \arrow[l] \\
        \Hom(\mid W \mid, \mid Y \mid) \arrow[u] &  \Hom(\mid Z \cup W \mid, \mid Y \mid) \arrow[l] \arrow[u]
    \end{tikzcd}
    \end{center}
    is homotopy cartesian, as desired.
    
    \item[II] The unique join-dimension 0 element in $\SCplx_V$ is the set containing the empty set, i.e., $\{\emptyset\} \in \SCplx_V$. We have $G(\{\emptyset\}) = \Hhom(\emptyset, |Y|) = \emptyset$, and $F(\{\emptyset\}) = | \Hhom(\{\emptyset\}, Y) | = |\{\emptyset\}| = \emptyset$.
    
    The join-dimension 1 elements in $\SCplx_V$ are the join-irreducible subcomplexes of $\Po(V)$. These are precisely the subcomplexes that are \emph{fully connected}, i.e., subcomplexes $Z$ of the form $Z = \Po(U)$ for any $U \subseteq V$ (we can think of this as the full simplex spanned by the vertices $U$).
    To see this, observe first that any subcomplex is a union of fully connected subcomlexes, so it suffices to show that fully connected subcomplexes are join-irreducible. Now, if $\Po(U) = A \cup B$ for some subcomplexes $A$ and $B$, then either $U \subseteq A$ or $U \subseteq B$, which implies that $\Po(U) \subseteq A$ or $\Po(U) \subseteq B$. Hence, $\Po(U)$ is join-irreducible.

    Firstly, as $|Z|$ is contractible when $Z$ is fully connected, we have
    \begin{equation*}
        G(Z) = \Hhom(|Z|, |Y|) \simeq \Hhom(*, |Y|) \cong |Y|
    \end{equation*}
    whenever $\jdim(Z) = 1$.

    We now show that also $F(Z) \simeq |Y|$ whenever $\jdim(Z) = 1$. Suppose $Z$ is fully connected, and write $Z \cong \Delta^n = \Po(\{0, \dots, n\})$. There is a simplicial map $H \colon \Delta^n \times \Delta^1 \to \Delta^n$ defined on vertices by
\begin{equation*}
    H(v, t) =
    \begin{cases}
        v, &\quad t = 0, \\ 
        0, &\quad t = 1. \\ 
    \end{cases}
\end{equation*}
    Applying $\Hhom(-,Y)$ gives a map
    \begin{equation*}
        H^* \colon \Hhom(\Delta^n, Y) \to \Hhom(\Delta^n \times \Delta^1, Y).
    \end{equation*}
    Applying the Hom-tensor adjunction two times, we get the map
    \begin{equation*}
        \Tilde H \colon \Hhom(\Delta^n, Y) \times \Delta^1 \to \Hhom(\Delta^n, Y) 
    \end{equation*}
    defined by $\Tilde H (f,t)(v) = H^*(f)(v,t) = f(H(v,t))$. Taking the geometric realization of $\Tilde H$ gives a deformation retract from $|\Hom(\Delta^n, Y)|$ to $|\Hom(*, Y)| \cong |Y|$, as desired.
\end{description}

\end{example}

\bibliographystyle{plain}
\bibliography{ref}

\end{document}